# PENALIZED LOG-LIKELIHOOD ESTIMATION FOR PARTLY LINEAR TRANSFORMATION MODELS WITH CURRENT STATUS DATA[1]

By Shuangge Ma and Michael R. Kosorok

*University of Wisconsin–Madison*

We consider partly linear transformation models applied to current status data. The unknown quantities are the transformation function, a linear regression parameter and a nonparametric regression effect. It is shown that the penalized MLE for the regression parameter is asymptotically normal and efficient and converges at the parametric rate, although the penalized MLE for the transformation function and nonparametric regression effect are only $n^{1/3}$ consistent. Inference for the regression parameter based on a block jackknife is investigated. We also study computational issues and demonstrate the proposed methodology with a simulation study. The transformation models and partly linear regression terms, coupled with new estimation and inference techniques, provide flexible alternatives to the Cox model for current status data analysis.

**1. Introduction.** Partly linear transformation models are flexible semiparametric regression models where a continuous outcome $U$, conditional on covariates $Z \in \mathbb{R}^d$ and $W \in \mathbb{R}$, is modeled as

$$(1.1) \qquad H(U) = \beta' Z + h(W) + e,$$

where $H$ is an unknown nondecreasing transformation, $h$ is an unknown smooth function, and $e$ has a known distribution $F$ with support $\mathbb{R}$. The setting we focus on in this paper is when $U$ is not observed directly, but only its current status is observed at a random censoring point $V \in \mathbb{R}$. More specifically, we observe $X = (V, \Delta, Z, W)$, where $\Delta = \mathbf{1}_{(U \leq V)}$. We assume that $U$ and $V$ are independent given $(Z, W)$. Although it is not difficult to

Received July 2003; revised October 2004.

[1]Supported in part by National Cancer Institute Grant CA75142.

*AMS 2000 subject classifications.* Primary 62G08, 60F05; secondary 62G20, 62B10.

*Key words and phrases.* Current status data, empirical processes, nonparametric regression, semiparametric efficiency, splines, transformation models.

extend our approach to allow multivariate $h$, we restrict our attention to univariate $h$ for simplicity of exposition.

Linear transformation models, having the form (1.1) without the nonparametric regression term $h$, have a long history. A parametric version of this model, with $H$ specified up to a finite-dimensional parameter vector, was initially investigated by Box and Cox [8]. Dabrowska and Doksum [13] studied the generalized odds-rate model for the two-sample problem, a special case of transformation models with nonparametric $H$. Cheng, Wei and Ying [10] proposed a class of estimating functions for the regression parameter $\beta$, with possibly right censored observations, and verified that the resulting estimator is asymptotically normal. Bickel and Ritov [7] developed efficient methods for linear transformation models with covariates in the uncensored setting.

The model (1.1) can be readily applied to a failure time $T$ by letting $U = \log T$. The partly linear Cox model, for example, is obtained by changing the sign of $\beta$ and $h$, letting $F(s) = 1 - \exp\{-e^s\}$, and by taking $H(u) = \log A(e^u)$, where $A$ is an unspecified integrated baseline hazard. More specifically, the hazard function of $T$, given the covariates $Z = z$ and $W = w$, is assumed to have the form

$$(1.2) \qquad \lambda(t|z,w) = \exp\{\tilde{\beta}'z + \tilde{h}(w)\}a(t),$$

where

$$\tilde{\beta} = -\beta, \qquad \tilde{h} = -h,$$

and where $a$ is the baseline hazard. This model has been studied for right censored data by Huang [23]. The partly linear proportional odds survival regression model has the same form only with $F(s) = e^s/(1 + e^s)$. If $Y$ is a current status time, or "case 1" interval censoring time, then letting $V = \log Y$ will yield the data structure described in the first paragraph.

Statistical methodology for current status data also has a long history. An important early example of current status data comes from tumor studies in animals, where the time of tumor onset is of interest, but not directly observable, as discussed in [15]. Current status data may occur due to study design or measurement limitations. Examples of such data arise in several fields, including demography, epidemiology, econometrics and bioassay. Research into statistical methods for current status data appears to have begun with the paper by Ayer, Brunk, Ewing, Reid and Silverman [2] on estimating a distribution function from a single sample. Other early approaches to current status data analysis include the use of generalized linear regression models [3, 14, 24, 50]. Andrews, van der Laan and Robins [1] investigate locally efficient estimators of parametric regression models with current status data.

Over the past decade, the fascinating asymptotic properties of estimators in the single sample case [19], and the similarity to more general kinds of interval censoring, have kindled significant interest in nonparametric approaches to estimation in current status regression [4, 16, 17, 20, 22, 28, 32, 34, 35, 36, 38, 40, 41, 42, 44, 45, 46, 52].

Our goal in this paper is to study estimation in the model (1.1), with special attention on inference for $\beta$. Our results extend the nonparametric likelihood-based approach of Huang [22, 23] in two ways. First, a smooth nonparametric covariate effect is added to the Cox model for current status data. Second, results are obtained for general transformation models with arbitrary but known residual distribution $F$. The approach we take is to use nonparametric maximum penalized log-likelihood.

Several interesting issues arise from carrying through this extension. First, the convergence rates for the estimators of $h$ and $H$, $\hat{h}_n$ and $\hat{H}_n$, appear to interfere with each other so as to require oversmoothing of $\hat{h}_n$ and thus force the convergence rates to both equal $O_p(n^{-1/3})$. Second, $\hat{H}_n$ has a bias which does not vanish asymptotically, even when no regression terms are present. This bias arises from an assumption we make that the support of the current status value $V$ is a finite interval. This assumption is also made by Huang [22]. In spite of this persistent bias, $\hat{H}_n$ is $L_2$ consistent and bounded in probability, and thus is sufficiently consistent to enable weak convergence of the estimator of $\beta$, $\hat{\beta}_n$. Third, inference for $\beta$ is challenging because estimation of the covariance directly is impractical and there exists no generally applicable method of inference for penalized estimators of partly linear models. The likelihood ratio expansion results of Murphy and van der Vaart [31] cannot be used in our setting since the penalized component of the objective function is larger than $O_p(n^{-1})$ and thus not negligible in the limit. To resolve the inference problem, we use a block jackknife estimator which is computationally simple and which applies, in general, to asymptotically linear statistics.

There is an interesting connection between the model (1.1) and semiparametric binary choice models studied in econometrics [9, 11, 12, 21, 26]. The model (1.1) can also be expressed as the probability of a consumer choosing "$\Delta = 1$" instead of "$\Delta = 0$," given covariates $(Z, W, V)$, via the expression

$$(1.3) \qquad P[\Delta = 1 | Z, W, V] = F(\beta' Z + h(W) + H(V)),$$

where $H$ is assumed to be a monotone covariate effect, $F$ is a known function, and the other covariate effects are as defined previously. Without the $H$ term, (1.3) is precisely the model studied by Härdle, Mammen and Müller [21]. This connection has also been observed in other settings of current status data study, for example, [1, 41, 42, 46].

The next section, Section 2, presents the data and model assumptions, along with several examples of residual distributions $F$ which satisfy the given requirements. The maximum penalized log-likelihood estimation procedure is presented in Section 3. Results on consistency of the estimators and the persistent bias in $\widehat{H}_n$ are given in Section 4. Section 5 presents results on rates of convergence for parameter estimators, and Section 6 presents asymptotic normality and efficiency of $\hat{\beta}_n$. A block jackknife inference procedure for $\beta$ is presented in Section 7. Several computational issues are discussed in Section 8, and a simulation study evaluating the moderate sample size performance of the proposed methods is given in Section 9. Proofs are given in Section 10.

## 2. The data setup and model assumptions.

2.1. *Data and model assumptions.* The data $\{X_i = (V_i, \Delta_i, Z_i, W_i), i = 1, \ldots, n\}$ consists of $n$ i.i.d. realizations of $X = (V, \Delta, Z, W)$, as described in the Introduction, generated by the model (1.1). Recall that $\Delta = \mathbf{1}_{(U \le V)}$, where $U$ is a real-valued outcome of interest. We make the following additional assumptions about the covariates and censoring distribution.

A1. (a) The covariate $Z$ belongs to a bounded subset $\mathcal{Z} \subset \mathbb{R}^d$. (b) The support for $W$ is $[a, b]$, where $-\infty < a < b < \infty$. (c) The support for $V$ is an interval $[l_v, u_v]$, where $-\infty < l_v < u_v < \infty$.

A2. $E \operatorname{var}[Z|V, W]$ is positive definite.

A3. With probability $> 0$, the conditional distribution of $W$ given $V$ dominates the unconditional distribution of $W$.

A4. $U$ and $V$ are independent given $(Z, W)$.

Define the function class $\Im_\nu \equiv \{h : [a, b] \mapsto \mathbb{R} \text{ with } J(h) < \infty\}$, where $J^2(h) \equiv \int_a^b [h^{(\nu)}(w)]^2 \, dw$ for some positive integer $\nu$, and where $h^{(j)}$ is the $j$th derivative of $h$. Also define the following subset of $\Im_\nu$ : $\Im_{\nu,0}^c \equiv \{h : h \in \Im_\nu, \sup_w |h(w)| < c \text{ and } E[h(W)] = 0\}$, where $c$ is a positive constant. When $c = \infty$, the superscript is omitted. We make the following additional model assumptions.

B1. The distribution of $U$ given the covariates has the transformation model form given in (1.1).

B2. The true regression parameter $\beta_0$ belongs to a known, bounded open subset $B_0$ of $\mathbb{R}^d$.

B3. The true nonparametric covariate effect $h_0 \in \Im_{\nu,0}^{c_0}$, for some known $c_0 < \infty$.

B4. The true transformation $H_0 : \mathbb{R} \mapsto \mathbb{R}$ is strictly monotone increasing and bounded on compacts.

B5. (a) The residual error distribution $F$ is known. (b) $F$ has first and second derivatives $f$ and $\dot{f}$, respectively, where the support of $f$ is $\mathbb{R}$ and where $|\dot{f}|$ is bounded. (c) For each compact $K \subset \mathbb{R}$, there exist constants $c \in (0, \infty)$ and $\alpha \in (1/2, 1]$ and an increasing isomorphic function $\xi : [0, 1] \mapsto [0, 1]$ so that

$$\sup_{s \in K} \sup_{u, v \in [0,1] \,:\, |\xi(u) - \xi(v)| \leq \epsilon} |F(F^{-1}(u) + s) - F(F^{-1}(v) + s)| \leq c\epsilon^{\alpha},$$

for all $\epsilon \in (0, 1)$. (d) $[f^2(v) - \dot{f}(v)F(v)] \wedge [f^2(v) + \dot{f}(v)(1 - F(v))] > 0$, for all $v \in \mathbb{R}$.

REMARK 1. In condition B3 it is assumed that the unknown nonparametric covariate effect $h_0$ is bounded by $c_0$. In the theoretical proofs and numerical calculations the exact value of $c_0$ is not necessary. Instead, only the boundedness condition is needed. The condition B5(c) is used to control the entropy of the model in order to obtain consistency. Condition B5(d) ensures convexity of the function $s \mapsto \delta \log F(s) + (1 - \delta) \log(1 - F(s))$ for $\delta = 0, 1$. This convexity is used in the proof of Theorem 2 below to establish that the estimator of the transformation function is bounded above and below in probability.

REMARK 2. It is not hard to verify that if $u \mapsto F(u)$ satisfies B5(b)–B5(d), then for any $a \in (0, \infty)$ and $b \in \mathbb{R}$, $u \mapsto F(au + b)$ also satisfies B5(b)–B5(d).

2.2. *Examples of transformation models.* The following are several examples of residual error distribution functions.

1. $F(u) = 1 - \exp[-e^u]$ is the extreme value distribution and corresponds to the complementary log-log transformation.
2. $F(u) = e^u[1 + e^u]^{-1}$ is the logistic distribution and corresponds to the logit transformation.
3. $F(u) = 1 - [1 + \gamma e^u]^{-1/\gamma}$ is a Pareto distribution with parameter $\gamma \in [0, \infty)$ and corresponds to the odds-rate transformation family. Taking the limit as $\gamma \downarrow 0$ yields the extreme value distribution, while $\gamma = 1$ gives the logistic distribution.
4. $F(u) = \Phi(u) \equiv (2\pi)^{-1/2} \int_{-\infty}^{u} \exp[-s^2/2] \, ds$ is the standard normal distribution which corresponds to the probit link.
5. $F(u) = \gamma[2\Gamma(1/\gamma)]^{-1} \int_{-\infty}^{u} \exp[-|s|^\gamma] \, ds$, for $\gamma \in [1, \infty)$, is a family of distributions which, after appropriate rescaling as justified in Remark 2 above, includes $\Phi(u)$ (corresponding to $\gamma = 2$).
6. $F(u) = 1/2 + \pi^{-1} \tan^{-1}(u)$ is the Cauchy distribution.

The following lemma gives a few examples which satisfy B5.

LEMMA 1. *Examples* 1–4 *and example* 5 *with* $\gamma \in (1, \infty)$ *satisfy conditions* B5(b)–B5(d).

REMARK 3.   Actually, many residual error distributions satisfy B5, but example 5 with $\gamma = 1$ and example 6 do not satisfy condition B5(d) [although they do satisfy B5(b) and B5(c)]. Consider first the Cauchy example. To see that B5(d) is not satisfied, note that, for sufficiently large $u$, $1 - F(u) = [\pi u]^{-1} + o(1/u)$. Since $f(u) = [\pi(1 + u^2)]^{-1}$ and $\dot{f}(u) = -2\pi u f^2(u)$, we have for $u > 0$ that $\dot{f}(u)(1 - F(u)) + f^2(u) = -f^2(u)(1 + o(1/u))$ and, thus, B5(d) is not satisfied. In example 5 with $\gamma = 1$, we have that $\dot{f}(u)(1 - F(u)) + f^2(u) = 0$ for all $u > 0$ and, thus, B5(d) is again not satisfied.

Condition B5(c) can sometimes be assessed directly for a given error distribution. In the proof of Lemma 1 this approach is taken to verify B5(c) for examples 1–3. In other settings the condition is hard to use directly, and the sufficient conditions given in the following lemma are more easily established. This approach is used to establish B5(c) for examples 4 and 5 [or $\gamma \in (1, \infty)$].

LEMMA 2.   Let $F$ satisfy condition B5(b) and:

(i) For every $\tau \in [0, \infty)$, there exists a $c_1^{(\tau)} \in [0, \infty)$ such that $\dot{f}(-u) < 0 < \dot{f}(u)$ and

$$f(-u)\dot{f}(-u + \tau) - \dot{f}(-u)f(-u + \tau) \leq 0 \leq f(u)\dot{f}(u - \tau) - \dot{f}(u)f(u - \tau),$$

for all $u \in (c_1^{(\tau)}, \infty)$.

(ii) For every $\tau \in [0, \infty)$, there exist a $c_2^{(\tau)} \in (0, \infty)$, an $\alpha_\tau \in (1/2, 1]$ and an increasing isomorphic function $\xi_*^{(\tau)} : [0, 1] \mapsto [0, 1]$ such that

$$F(F^{-1}(\xi_*^{(\tau)}\{\epsilon\}) + \tau) \vee [1 - F(F^{-1}(\xi_*^{(\tau)}\{1 - \epsilon\}) - \tau)] \leq c_1^{(\tau)}\epsilon^{\alpha_\tau},$$

for all $\epsilon$ sufficiently small.

Then $F$ satisfies B5(c).

3.  **Maximum penalized log-likelihood estimation.**  Under model (1.1) the log-likelihood for a single observation at $X = x \equiv (v, \delta, z, w)$ for the parameter choice $(\beta, h, H)$ is

$$
\begin{aligned}
(3.1) \quad l(x; \beta, h, H) &\equiv \delta \log\{F[\beta'z + h(w) + H(v)]\} \\
&\quad + (1 - \delta) \log\{1 - F[\beta'z + h(w) + H(v)]\}.
\end{aligned}
$$

Intuitively, we will need some mechanism to control the smoothness of estimates of $h_0$. One approach is to use sieve estimates with assumptions on the derivatives, as in [23]. However, we use instead a penalized approach based on splines. An advantage is that the degree of smoothness is controlled by a single number, the penalty term. Specifically, we use the following penalized log-likelihood function based on $n$ observations:

$$l_n^p(\beta, h, H) \equiv \mathbb{P}_n l(x; \beta, h, H) - \lambda_n^2 J^2(h),$$

where $\mathbb{P}_n$ is the empirical measure based on the $n$ observations $X_1, \ldots, X_n$ and $\lambda_n$ is the (possibly data-dependent) smoothing parameter. We maximize $l_n^p$ under the constraints that $\beta \in \bar{B}_0$, $h \in \Im_\nu$ with $|h| \leq c_0$ and $\mathbb{P}_n h(W) = 0$, and that $H$ is nondecreasing, where $\bar{B}$ denotes the closure of the set $B$. We also define $l_n(\beta, h, H) \equiv \mathbb{P}_n l(x; \beta, h, H)$ to be the ordinary log-likelihood function. The maximum penalized log-likelihood estimators $\hat{\beta}_n$, $\hat{h}_n$ and $\widehat{H}_n$ are the ones that maximize the penalized log-likelihood function under the stated constraints. We make the following assumption about the penalty term $\lambda_n$:

C. $\lambda_n = O_p(n^{-1/3})$ and $\lambda_n^{-1} = O_p(n^{1/3})$.

REMARK 4. The tuning parameter $\lambda_n$ is usually selected through certain cross validation techniques (see [49] for reference). However, for the asymptotics to hold, we only require $\lambda_n$ to be of the correct order. One way of achieving this is to simply set $\lambda_n = n^{-1/3}$, as we do for the simulation studies shown in Section 9. In addition, while our theoretical results require specification of $c_0$, it appears from our experience that $c_0$ does not actually need to be specified for implementation with moderate sample sizes. Numerical studies show for finite sample sizes that this strategy usually yields satisfactory results. Note that we are oversmoothing in our choice of penalty term. This is a consequence of interference from estimation of $H$. In Theorem 4 presented in Section 5, and in the proof of the theorem given in Section 10, it is demonstrated that the entropy for the model is driven by the entropy of the class of monotone functions which parameterizes $H$, resulting in an overall rate of $O_p(n^{-1/3})$. Although we can refine the rate for estimation of $\beta$ to $O_p(n^{-1/2})$, it is unclear how to improve the rate of $O_p(n^{-1/3})$ for estimating $h$ when $\nu > 1$. It is an open question whether achieving the optimal rate for $h$ in model (1.1) is even possible using penalized log-likelihood.

REMARK 5. Let $y_{(1)}, \ldots, y_{(n)}$ be the order statistics of $y_1, \ldots, y_n$. Let $\delta_{(i)}$, $w_{(i)}$ and $z_{(i)}$ correspond to $y_{(i)}$. Since only the values of $H$ at $y_{(i)}$ matter in the log-likelihood function, we will take the maximum likelihood estimator $\widehat{H}_n$ as the right-continuous nondecreasing step functions with jump points at $y_{(i)}$.

REMARK 6. From time to time, it will be convenient to assume that $\delta_{(1)} = 1$ and $\delta_{(n)} = 0$. Such an assumption usually has little impact on results. To see this, note that if $\delta_{(1)} = 0$, then this observation makes a contribution of zero to the log-likelihood function after maximizing over $H$. Similarly, if $\delta_{(n)} = 1$, then the corresponding observation also makes no contribution to the log-likelihood. We will make our use of this assumption clear. Further discussion about this assumption can be found in Section 3.1 of [19].

**4. Consistency.** The following lemma establishes existence of the maximum penalized log-likelihood estimators.

LEMMA 3. *Under the assumptions* A1–A4, B1–B5, *the assumption of Remark* 6, *and provided* $0 < \lambda_n < \infty$, *a maximum penalized log-likelihood estimator* $\hat{\psi}_n \equiv (\hat{\beta}_n, \hat{h}_n, \widehat{H}_n)$ *exists, with* $\hat{\beta}_n \in \bar{B}_0$, $\hat{h}_n \in \Im_\nu$, $\sup_{s \in [a,b]} |\hat{h}_n(s)| \leq c_0$, $\mathbb{P}_n \hat{h}_n(W) = 0$ *and* $-\infty < \widehat{H}_n(y_{(1)}) \leq \widehat{H}_n(y_{(n)}) < \infty$.

REMARK 7. Provided $\delta_{(j)} = 1$ and $\delta_{(j+1)} = 0$ for some $j \in \{1, \ldots, n-1\}$, Lemma 3 implies that $-\infty < \widehat{H}_n(y_{(k)}) \leq \widehat{H}_n(y_{(l)}) < \infty$, where $k = \inf\{j \in \{1, \ldots, n\} : \delta_{(j)} = 1\}$ and $l = \sup\{j \in \{1, \ldots, n\} : \delta_{(j)} = 0\}$. Note that the log-likelihood contributions for the observations corresponding to $\delta_{(j)}$ for $j \notin \{k, \ldots, l\}$ are zero.

Define the following distance between parameters $(h_1, H_1)$ and $(h_2, H_2)$: $d_F((h_1, H_1), (h_2, H_2)) \equiv \|h_1 - h_2\|_2 + \|H_1 - H_2\|_{F,2}$, where $\|h_1 - h_2\|_2 \equiv [\int_a^b |h_1(w) - h_2(w)|^2 \, dw]^{1/2}$ and $\|H_1 - H_2\|_{F,2} \equiv [\int_{l_v}^{u_v} |F(H_1(v)) - F(H_2(v))|^2 \, dv]^{1/2}$. Note that since $F$ has a bounded derivative, $\| \cdot \|_{F,2}$ is somewhat weaker than the usual $L_2$ norm. We use $\| \cdot \|$ to denote Euclidean distance. The following is the main consistency result.

THEOREM 1. *Under the assumptions* A1–A4, B1–B5 *and* C, $d_F((\hat{h}_n, \widehat{H}_n), (h_0, H_0)) + \|\hat{\beta}_n - \beta_0\| = o_p(1)$.

REMARK 8. We will show in Section 5 that $J(\hat{h}_n) = O_p(1)$ and, thus, Theorem 1 combined with condition A1(b) will imply $\sup_{w \in [a,b]} |\hat{h}_n(w) - h_0(w)| = o_p(1)$, since the $\hat{h}_n$'s are smooth functions defined on a compact set with asymptotically bounded first-order derivatives.

Under the provision given in Remark 7, and with $k$ and $l$ as given in that remark, define

$$\widetilde{H}_n(t) \equiv \begin{cases} \widehat{H}_n(y_{(k)}), & t \in [l_v, y_{(k)}), \\ \widehat{H}_n(t), & t \in [y_{(k)}, y_{(l)}], \\ \widehat{H}_n(y_{(l)}), & \text{otherwise}, \end{cases}$$

and let $\widetilde{H}_n = 0$ if the provision is not met. Theorem 1 clearly implies that $\|\widetilde{H}_n - H_0\|_{F,2} = o_p(1)$; and although $L_2$-convergence does not imply uniform convergence, the following theorem ensures that $|\widetilde{H}_n|$ is bounded.

THEOREM 2. *Under the conditions of Theorem 1, we have* $|\widetilde{H}_n(l_v)| \vee |\widetilde{H}_n(u_v)| = O_p(1)$.

REMARK 9. Theorems 1 and 2 jointly imply that $\int_{l_v}^{u_v} |\widetilde{H}_n(v) - H_0(v)|^2 \, dv = o_p(1)$. It appears that this cannot be strengthened to uniform consistency. In the following theorem, we prove in the setting where covariate effects are not included in the model that uniform consistency of $\widetilde{H}_n$ is impossible under conditions A1(c) and B4 and when the distribution of $V$ is continuous.

THEOREM 3. *Assume the nonparametric current status model* [*model* (1.1) *without $\beta$ or $h$*] *holds, and let $G_0$ denote the unknown true distribution of $U$. Assume condition* A1(c) *holds for the current status time and that $G_0$ is continuous with $0 < G_0(l_v) \leq G_0(u_v) < 1$* (*essentially condition* B4). *Assume also that the distribution of $V$ is continuous. Let $\widehat{G}_n$ be the nonparametric maximum likelihood estimator of $G_0$* (*no penalization is involved since $h$ is not in the model*), *and assume that the condition in Remark* 6 *holds. Then:*

(i) *for any $\epsilon \in (0, G_0(l_v))$, $\liminf_{n \to \infty} P\{\widehat{G}_n(V_{(1)}) \leq G_0(l_v) - \epsilon\} > 0$, and*

(ii) *for any $\epsilon \in (0, \{1 - G_0(u_v)\})$, $\liminf_{n \to \infty} P\{\widehat{G}_n(V_{(n)}) \geq G_0(u_v) + \epsilon\} > 0$.*

REMARK 10. Through simulation studies involving up to 20,000 observations, we have verified that the bias predicted in Theorem 3 does persist as $n$ gets large. This does not contradict the uniform consistency result in Section II.4.1 of [19]. In their consistency proof, they require the probability measure for the current status value $V$ to dominate the distribution of $U$. Thus, the total variation of $\widehat{G}_n$ over the support of $V$ is bounded by 1, and the total variation of $G_0$ over the support of $V$ is equal to 1. Thus, in this setting, $L_2$-convergence of $\widehat{G}_n$ with respect to the measure $G_0$ will indeed imply uniform convergence.

REMARK 11. By expression (10.5) in the proof of Theorem 1, and by the results of Theorems 1 and 2, it is clear that if

$$(4.1) \qquad P(V = l_v) \wedge P(V = u_v) > 0,$$

then $\widetilde{H}_n$ is uniformly consistent for $H_0$ over the interval $[l_v, u_v]$. Unfortunately, the assumption (4.1) does not seem to be very realistic. Fortunately, for the results that follow, $L_2$ convergence of $\widetilde{H}_n$ is sufficient.

**5. Rates of convergence.** For this section we will use the usual $L_2$ distance between parameters $(h_1, H_1)$ and $(h_2, H_2)$: $d((h_1, H_1), (h_2, H_2)) \equiv \|h_1 - h_2\|_2 + \|H_1 - H_2\|_2$, where $\|H_1 - H_2\|_2 \equiv [\int_{l_v}^{u_v} |H_1(v) - H_2(v)|^2 \, dv]^{1/2}$. We now establish the rate of convergence for all parameters.

THEOREM 4. *Rate of convergence: Suppose that assumptions* A1–A5, B1–B4 *and* C *are satisfied. Then* $J(\hat{h}_n) = O_p(1)$ *and* $\|\hat{\beta}_n - \beta_0\| + d((\hat{h}_n, \widetilde{H}_n), (h_0, H_0)) = O_p(n^{-1/3})$.

REMARK 12. In ordinary spline settings, $\lambda_n = O_p(n^{-\nu/(2\nu+1)})$ and the optimal rate of convergence for a smooth function of the covariate $w$ is $O_p(\lambda_n)$ (see, e.g., [49]). Typically, $\lambda_n$ is a data driven smoothing parameter selected by cross validation, GCV, $L_P$ or $L_C$, as discussed in [27]. However, in our case the rate of convergence for $\widehat{H}_n$ cannot exceed $n^{1/3}$, which slows down the convergence of $\hat{h}_n$. In particular, the rate for $\hat{h}_n$ does not achieve the optimal convergence rate of $O_p(\lambda_n)$ when $\nu > 1$, as is achieved in [30]. It is also worth pointing out that it appears we cannot achieve a better convergence rate by modifying the smoothness assumptions.

REMARK 13. As discussed in Section 1, a special case of the general transformation models is the proportional hazard model. It is shown by Groeneboom and Wellner [19] that the convergence rate of the NPMLE of a distribution function is $n^{1/3}$ for current status data. So we have shown that, under reasonable model assumptions, the convergence rate of $\widetilde{H}_n$ achieves the optimal rate, despite the presence of an additional infinite-dimensional parameter $h$.

In some cases it is reasonable to assume that the transformation function is also continuously differentiable. In [29] penalized estimation of the transformation function is investigated under certain smoothness assumptions. A sharper convergence rate can be achieved if we assume the transformation function belongs to a certain Sobolev space. In our case, if we assume the transformation function $H_0$ and the nonparametric effect $h_0$ belong to the same Sobolev space, then we can achieve optimal convergence rates for both parameters by using doubly penalized estimators. However, if it is assumed that $h$ and $H$ belong to different Sobolev spaces, then it is unclear whether we can achieve the optimal convergence rate for both estimators.

## 6. Weak convergence of the parametric covariate effect.

6.1. *Information calculation.* It is well known that in most parametric models we can estimate the finite-dimensional parameter at the $n^{1/2}$ convergence rate. However, this is not necessarily true for semiparametric models. A necessary condition is that we have positive information, which is not a trivial condition. Next we show that, for our model, indeed, we have positive information. The following extra model assumptions will be needed:

D1. $\exists \tilde{h} \in \Im_{\nu,0}$ such that

$$E\left\{ [Z - \tilde{h}(W)] \times \left[ h(W) - \frac{E(h(W)Q_{\psi_0}^2(X)|V)}{E(Q_{\psi_0}^2(X)|V)} \right] Q_{\psi_0}^2(X) \right\} = 0$$

for every $h \in \Im_{\nu,0}$, where

$$Q_\psi(x) \equiv f(v)\left\{ \frac{\delta}{F(\theta_\psi(x))} - \frac{1-\delta}{1-F(\theta_\psi(x))} \right\},$$

$\psi \equiv (\beta, h, H), \psi_0 \equiv (\beta_0, h_0, H_0), \theta_\psi(x) \equiv \beta' z + h(w) + H(v)$ and

$$\tilde{q}(v) \equiv \frac{E(ZQ_{\psi_0}^2(X)|V=v)}{E(Q_{\psi_0}^2(X)|V=v)} - \frac{E(\tilde{h}(W)Q_{\psi_0}^2(X)|V=v)}{E(Q_{\psi_0}^2(X)|V=v)}$$

has a derivative which is uniformly bounded on $[l_v, u_v]$.

D2. $I_0 \equiv E(\dot{l}\dot{l}')$ is positive definite, where $\dot{l} \equiv \{Z - \tilde{h}(W) - \tilde{q}(V)\}Q_{\psi_0}(X)$.

Before giving the main result of this section, Theorem 5, we present a lemma which provides sufficient conditions for achieving D1 and determining $\tilde{h}$. Let $D_0(v,w) \equiv E[Q_{\psi_0}^2(X)|V=v, W=w], D_1(v,w) \equiv E[ZQ_{\psi_0}^2(X)|V=v, W=w], D_{01}(v) \equiv E[Q_{\psi_0}^2(X)|V=v]$ and $D_{02}(w) \equiv E[Q_{\psi_0}^2(X)|W=w]$; and define $\mathcal{S}_0$ to be the class of functions $g : [l_v, u_v] \times [a, b] \mapsto \mathbb{R}$ with $E[g^2(V,W) \times D_0(V,W)] < \infty$. For any $g \in \mathcal{S}_0$, let $\Pi_1$ be the projection operator

$$g \mapsto \frac{E[g(V,W)Q_{\psi_0}^2(X)|V=v]}{E[Q_{\psi_0}^2(X)|V=v]}$$

and let $\Pi_2$ be the projection operator

$$g \mapsto \frac{E[g(V,W)Q_{\psi_0}^2(X)|W=w]}{E[Q_{\psi_0}^2(X)|W=w]}.$$

Also define $D^* \equiv D_1/D_0, D_1^* \equiv \Pi_1 D^*, D_2^* \equiv \Pi_2 D^*, R(v,w) \equiv D_0(v,w)/D_{01}(v), S(v,w) \equiv D_0(v,w)/D_{02}(w), \dot{D}_1^*(v) \equiv (\partial/(\partial v))D_1^*(v), \dot{R}_1(v,w) \equiv (\partial/(\partial v))R(v,w), D_2^{*(\nu)}(w) \equiv (\partial/(\partial w))^\nu D_2^*(w)$ and $S_2^{(\nu)}(v,w) \equiv (\partial/(\partial w))^\nu S(v,w)$.

LEMMA 4. *Assume the model conditions of Section* 2.1. *Also assume that $V$ and $W$ are independent with the density of $W$ being bounded above and below on $[a, b]$, that $\dot{D}_1^*(v)$ and $E[\dot{R}_1(v, W)]^2$ are uniformly bounded on $[l_v, u_v]$, and that both $E[D_2^{*(\nu)}(W)]^2 < \infty$ and $E[S_2^{*(\nu)}(V, W)]^2 < \infty$. Then* D1 *is satisfied with $\tilde{h} = h^* - Eh^*(W)$, where $h^* \equiv \lim_{m \to \infty} [\sum_{j=0}^{m} (\Pi_2 \Pi_1)^j] \Pi_2 (1 - \Pi_1) D^*$.*

REMARK 14. We note that while the conditions of Lemma 4 are somewhat stronger than our previous model assumptions, the conditions are still reasonable. The most arduous of the new assumptions involve bounding the derivatives of several well-defined functions. Note that the denominators of the ratios that define these functions are $D_{01}(v)$ and $D_{02}(w)$. Since $D_{01}$ and $D_{02}$ are bounded below on $[l_v, u_v]$ and $[a, b]$, respectively, the task of bounding the necessary derivatives is simplified somewhat.

THEOREM 5. *Calculation of efficient information*: *Under model assumptions* A1–A4, B1–B5 *and* D1–D2, $I_0$ *is the efficient information matrix for $\beta$.*

REMARK 15. Knowledge about the degree of smoothness for $h$ cannot currently be utilized to improve the rate of convergence for $\hat{h}_n$ or the asymptotic precision of $\hat{\beta}_n$. This is a consequence of the fact that the Sobolev space $\Im_{\nu_1}$ is dense in $\Im_{\nu_2}$ when $\nu_1 > \nu_2$. However, it is unclear whether such knowledge can result in small sample improvements in the accuracy of $\hat{\beta}_n$.

### 6.2. *Asymptotic normality and efficiency.*

THEOREM 6. *Asymptotic normality and efficiency*: *Assume conditions* A1–A4, B1–B5, C *and* D1–D2. *Also assume that the inverse of $H_0$, $H_0^{-1}$, has a derivative bounded on compacts. Then $n^{1/2}(\hat{\beta}_n - \beta_0) = I_0^{-1} \sqrt{n} \, \mathbb{P}_n \tilde{l} + o_p(1) \overset{d}{\to} N(0, I_0^{-1})$.*

Since $\hat{\beta}_n$ is asymptotically linear with efficient influence function, and the model is sufficiently smooth (Hellinger differentiable), it is asymptotically efficient in the sense that any regular estimator has asymptotic covariance matrix no less than that of $\hat{\beta}_n$ [6]. The additional assumption on the smoothness of $H_0^{-1}$ is needed to construct an approximately least-favorable submodel (see Section 25.11 of [47]) under which the given estimator satisfies the efficient score equation.

Another important issue for statistical modeling with current status data is the degree of robustness achieved under model misspecification. Yu and van der Laan [53] investigate doubly robust estimation in longitudinal marginal structural models. Their results allow one to construct locally efficient estimators of the regression parameters, under the misspecification of either the unknown regression function or the conditional distribution of the linear variables. The question of whether double robustness holds in our setting is also of interest, but is beyond the scope of the current paper.

**7. The block jackknife.** One potential method of inference for $\beta$ is to use the nonparametric bootstrap. Unfortunately, there is no sufficiently general theory, as far as we are aware, available for the nonparametric bootstrap in the penalized maximum likelihood estimation setting. Alternative approaches are the $m$ within $n$ bootstrap (see [5]) or subsampling (see [33]). Since $\sqrt{n}(\hat{\beta}_n - \beta_0)$ has a continuous limiting distribution $\mathcal{L}$, Theorem 2.1 of [33] yields that the $m$ out of $n$ subsampling bootstrap converges—conditionally on the data—to the same distribution $\mathcal{L}$, provided $m/n \to 0$ and $m \to \infty$ as $n \to \infty$. Because of the requirement that $m \to \infty$ as $n \to \infty$, the subsampling bootstrap potentially involves many calculations of the estimator. Fortunately, the asymptotic linearity given in Theorem 6 can be used to formulate a computationally simpler method of inference.

Let $\tilde{\beta}_n$ be any asymptotically linear estimator of a parameter $\beta_0 \in \mathbb{R}^d$, based on an i.i.d. sample $X_1, \ldots, X_n$, having square-integrable influence function $\phi$ for which $E[\phi\phi^T]$ is nonsingular. Let $m$ be a fixed integer $> d$, and, for each $n \geq m$, define $k_{m,n}$ to be the largest integer satisfying $mk_{m,n} \leq n$. Also define $N_{m,n} \equiv mk_{m,n}$. For the data $X_1, \ldots, X_n$, compute $\tilde{\beta}_n$ based on the proposed estimation method and randomly sample $N_{m,n}$ out of the $n$ observations without replacement, to obtain $X_1^*, \ldots, X_{N_{m,n}}^*$. For $j = 1, \ldots, m$, let $\tilde{\beta}_{n,j}^*$ be the estimate of $\beta$ based on the observations $X_1^*, \ldots, X_{N_{m,n}}^*$ after omitting $X_j^*, X_{m+j}^*, X_{2m+j}^*, \ldots, X_{(k_{m,n}-1)m+j}^*$. Compute

$$\bar{\beta}_n^* \equiv m^{-1} \sum_{j=1}^m \tilde{\beta}_{n,j}^* \quad \text{and} \quad S_n^* \equiv (m-1)k_{m,n} \sum_{j=1}^m (\tilde{\beta}_{n,j}^* - \bar{\beta}_n^*)(\tilde{\beta}_{n,j}^* - \bar{\beta}_n^*)^T.$$

The following lemma provides a method of obtaining asymptotically valid confidence ellipses for $\beta_0$.

LEMMA 5. *Let $\tilde{\beta}_n$ be an estimator of $\beta_0 \in \mathbb{R}^d$, based on an i.i.d. sample $X_1, \ldots, X_n$, which satisfies $n^{1/2}(\tilde{\beta}_n - \beta_0) = \sqrt{n}\mathbb{P}_n\phi + o_p(1)$, where $E[\phi\phi^T]$ is nonsingular. Then $n(\tilde{\beta}_n - \beta_0)^T[S_n^*]^{-1}(\tilde{\beta}_n - \beta_0)$ converges weakly to $d(m-1) \times F_{d,m-d}/(m-d)$, where $F_{r,s}$ has an F distribution with degrees of freedom $r$ and $s$.*

The key to the proof of Lemma 5 is the simultaneous validity of the asymptotic linearity expansion for all of the jackknife estimates. The details of the proof are given in the Appendix.

REMARK 16. Let $S_n^{**}$ be $S_n^*$ with $\bar{\beta}_p^*$ replaced by the estimator of $\beta$ based on $X_1^*, \ldots, X_{N_{m,n}}^*$ (which we denote $\tilde{\beta}_n^*$). Using arguments in the proof of Lemma 5, it is straightforward to show that replacing $S_n^*$ with $S_n^{**}$ will not affect the conclusions of Lemma 5. Those same arguments also lead to the conclusion that one cannot, in general, replace $\bar{\beta}_n^*$ with $\hat{\beta}_n$, except when $N_{m,n} = n$ (in which case $\hat{\beta}_n = \tilde{\beta}_n^*$).

REMARK 17. The block jackknife procedure only requires computing the estimator $m + 1$ times (or $m + 2$ times when $S_n^{**}$ is used as discussed in Remark 16). In our simulation studies we have found that for the proposed estimator $m = 10$ works for sample sizes $n = 400$ and 1600. The fact that $m$ remains fixed as $n \to \infty$ in the proposed approach results in a potentially significant computational savings over subsampling, which requires $m$ to grow increasingly large as $n \to \infty$. A potential challenge for the proposed approach is in choosing $m$ for a given data set. The larger $m$ is, the larger the denominator degrees of freedom in $F_{d,m-d}$ and the tighter the confidence ellipsoid. On the other hand, $m$ cannot be so large that the asymptotic linearity of Theorem 6 does not hold simultaneously for all jackknife components.

## 8. Computational techniques.

8.1. *Overall strategy.* Computationally, finding the penalized MLE for $\beta$, $h$ and $H$ is a maximization problem subject to the boundedness constraint for $h$ and the nondecreasing constraint for $H$. It is unlikely that there exists an analytic solution for this model. So we propose the following iterative maximization technique.

S1. For a given $\hat{\beta}_n^{(k)}$ and $\hat{h}_n^{(k)}$, estimate $\hat{H}_n^{(k)}$ by maximizing $l_n^p$ with respect to $H$ under the constraint that $\hat{H}_n^{(k)}$ is a nondecreasing step function.

S2. For a given $\hat{H}_n^{(k)}$, find $\hat{\beta}_n^{(k+1)}$ and $\hat{h}_n^{(k+1)}$ that maximize the penalized log-likelihood function.

We have found that almost any reasonable initial values will work. The above two-step maximization procedure is repeated until certain convergence criteria are satisfied. The global convexity of the log-likelihood function guarantees that we can reach the maximum by the above technique. For step S2, our experience indicates that, in applications involving moderate sample sizes, specification of $c_0$ is not needed and $\lambda_n = n^{-1/3}$ appears to work most of the time. Perhaps using cross validation to choose $\lambda_n$ may improve the performance of the estimator in some settings, but evaluating this issue requires further study and is beyond the scope of the current paper.

REMARK 18. After finite iterations, what we get is not exactly the penalized MLE. However, a very nice property of the efficiency theorem is that we only need approximate maximization to achieve asymptotic efficiency for $\beta$, that is, $\mathbb{P}_n \tilde{k}_{\hat{\beta}_n, \hat{h}_n, \widehat{H}_n} = o_p(n^{-1/2})$, where $\tilde{k}$ is the estimating function as defined in the proof of Theorem 6.

8.2. *Sieve approximation for the nonparametric covariate effect.* For the special case of $\nu = 2$, we can use a cubic spline for estimating $h$. Suppose a function $\hat{h}_n^*$ maximizes the penalized log-likelihood function. Then there exists a cubic spline function $\hat{h}_n$ such that $\hat{h}_n^*(w_i) = \hat{h}_n(w_i)$ for $i = 1, \ldots, n$ and $J(\hat{h}_n) \le J(\hat{h}_n^*)$. For a proof, see page 18 of [18]. The number of basis functions of a cubic spline increases at the rate $O(n)$. Computationally this can be quite time consuming for a moderate or large set of observations. Hence, we take a computational sieve approach suggested by Xiang and Wahba [51], which states that an estimate with the number of basis functions growing at least at the rate $n^{1/5}$ can achieve the same asymptotic precision as the full space. The $K$-mean clustering technique (see [25] for reference) is used to select the proper positions of knots, and $B$-spline basis functions are utilized. Because of the accuracy of this sieve approximation and the fact that the degree of smoothness is still controlled by the penalty term, the theoretical properties of the resulting estimators should be unmodified from the previously derived theory.

8.3. *Estimation of the transformation function.* The maximization over the nondecreasing function $H$ can be solved by some commercial software package, such as NPSOL. However, we show that the cumulative sum diagram approach, as discussed by Groeneboom and Wellner [19] and Huang [22], works for general transformation models. First we observe the following properties of $\widehat{H}_n$.

LEMMA 6. *Assume that $\delta_{(1)} = 1$, $\delta_{(n)} = 0$. Then for any fixed $\hat{\beta}_n$ and $\hat{h}_n$ the maximum likelihood estimator $\widehat{H}_n$ satisfies*

$$
\begin{aligned}
(8.1) \quad \sum_{j \ge i} \Bigg\{ & \delta_i \frac{f(\hat{\beta}_n' z_{(j)} + \hat{h}_n(w_{(j)}) + \widehat{H}_n(v_{(j)}))}{F(\hat{\beta}_n' z_{(j)} + \hat{h}_n(w_{(j)}) + \widehat{H}_n(v_{(j)}))} \\
& - (1 - \delta_i) \frac{f(\hat{\beta}_n' z_{(j)} + \hat{h}_n(w_{(j)}) + \widehat{H}_n(v_{(j)}))}{1 - F(\hat{\beta}_n' z_{(j)} + \hat{h}_n(w_{(j)}) + \widehat{H}_n(v_{(j)}))} \Bigg\} \le 0
\end{aligned}
$$

*for $i = 1, \ldots, n$, and*

$$
\sum_{i=1}^n \Bigg\{ \delta_i \frac{f(\hat{\beta}_n' z_{(i)} + \hat{h}_n(w_{(i)}) + \widehat{H}_n(v_{(i)}))}{F(\hat{\beta}_n' z_{(i)} + \hat{h}_n(w_{(i)}) + \widehat{H}_n(v_{(i)}))}
$$

(8.2)
$$-(1-\delta_i)\frac{f(\hat{\beta}_n' z_{(i)} + \hat{h}_n(w_{(i)}) + \widehat{H}_n(v_{(i)}))}{1-F(\hat{\beta}_n' z_{(i)} + \hat{h}_n(w_{(i)}) + \widehat{H}_n(v_{(i)}))}\bigg\}\widehat{H}_n(v_{(i)}) = 0.$$

This lemma can be proved in a manner similar to Proposition 1.1 of [19], and we omit the details. These properties motivate us to consider the following iterative, but computationally efficient algorithm.

Define the process $W_H$, $G_H$ and $O_H$ by

$$W_H(v) = \int_{v' \in [0,v]} \left\{\Delta\frac{f(\theta_\psi)}{F(\theta_\psi)} - (1-\Delta)\frac{f(\theta_\psi)}{1-F(\theta_\psi)}\right\} dR_n,$$

$$G_H(v) = \int_{v' \in [0,v]} \Delta\frac{f^2(\theta_\psi)}{F^2(\theta_\psi)(1-F(\theta_\psi))} dR_n,$$

$$O_H(v) = W_H(v) + \int H(v') dG_H,$$

where $v \geq 0$ and $R_n$ is the unobserved empirical measure of $(U, V, Z, W)$.

LEMMA 7. *Self-induced calculation of* $\widehat{H}_n$. *Assume that* $\delta_{(1)} = 1$ *and* $\delta_{(n)} = 0$. *Then for any fixed* $\hat{\beta}_n$ *and* $\hat{h}_n$, *the maximum likelihood estimator* $\widehat{H}_n$ *is the left derivative of the greatest convex minorant of the "self-induced" cumulative sum diagram, consisting of the points* $(G_{\widehat{H}_n}(V_{(j)}), O_{\widehat{H}_n}(V_{(j)}))$ *and the origin* $(0,0)$.

The proof of Lemma 7 is analogous to that of Proposition 1.4 of [19] and is omitted. This lemma gives an iterative procedure for finding $\widehat{H}_n$, as discussed in [19]. Suppose $\widehat{H}_n^{(k)}$ is the result of the $k$th iteration; then $\widehat{H}_n^{(k+1)}$ is computed as the left derivative of the greatest convex minorant of the cumulative sum diagram, consisting of the points $(G_{\widehat{H}_n^{(k)}}(V_{(j)}), O_{\widehat{H}_n^{(k)}}(V_{(j)}))$ and the origin $(0, 0)$.

**9. Simulation study.** To evaluate the finite-sample performance of our estimators, we conduct a small simulation study with current status data for the partly linear Cox model. As discussed in Section 1, the Cox model is a special case of general transformation models, where $H(u) = \log(A(e^u))$ and $F(s) = 1 - \exp(-e^s)$. The event times are generated from equation (1.2), with regression coefficients $\beta_1 = 0.3$ and $\beta_2 = 0.25$. The covariate $Z_1$ is Uniform$[0.5, 1.5]$ and $Z_2$ is Bernoulli with probability of success 0.5. For simplicity we take $h(w) = \sin(w/1.2 - 1) - k_0$, with $W$ Uniform$[1, 10]$ and $k_0 = 0.06516$, and $A(u) = e^{k_0}(\exp(u/3) - 1)$. Censoring times are standard exponentially distributed conditional on being in the interval $[0.2, 1.8]$. For computational simplicity we leave $c_0$ unspecified and do not use a data driven mechanism to select $\lambda_n$. Instead, we use $\lambda_n = n^{-1/3}$. $\hat{h}_n$ and $\hat{A}_n$ are estimated based on the computational strategies discussed in Section 8. We simulate 200 realizations for sample sizes equal to 400 and 1600.

The summary statistics for our estimators are shown in Table 1. It can be seen that the sample means are quite close to the true values. The sample standard deviation for $\hat{\beta}_1$ based on sample size 400 is 0.284, compared with 0.139 for sample size 1600, resulting in a ratio of 2.04. According to the asymptotic normality property (Theorem 6), the ratio should be 2. Hence, the ratio estimated from the simulations matches the theory quite well in this instance. The same property can be observed for estimators of $\hat{\beta}_2$. Inference based on the block jackknife, as discussed in Section 7, with the modification given in Remark 16, is also presented in this table. The 95% confidence intervals generally have coverage within two Monte Carlo standard errors $(0.03 = 2\sqrt{0.05 \times 0.95/200})$ of 0.95, except when $m = 40$ and the sample size $n = 400$. This is possibly because $m$ is too large for the asymptotic linearity property to hold simultaneously for all $m$ block jackknifes at this sample size, as hinted at in Remark 17.

Histograms of $\hat{\beta}_1$ and $\hat{\beta}_2$ and a plot of $\hat{\beta}_1$ versus $\hat{\beta}_2$ are shown in Figure 1. We can see clearly that the marginal distributions and the joint distribution of $\hat{\beta}_1$ and $\hat{\beta}_2$ appear to be Gaussian. Estimates and pointwise 95% confidence intervals for $h$ and $A$ based on sample size 1600 and 200 realizations are shown in Figure 2 and Figure 3. It can be seen that true values for $h$ and $H$ both lie in the 95% pointwise confidence intervals.

## 10. Proofs.

TABLE 1

*Simulation results for the partly additive Cox model with current status data. Sample sizes are equal to 400 and 1600. Sample means, standard deviations and confidence region coverages are based on 200 replicates. Confidence intervals are based on the block jackknife with $m = 10$ and 40 blocks. The true values of the regression parameters are $\beta_1 = 0.3$ and $\beta_2 = 0.25$*

|  |  | Sample size 400 | Sample size 1600 |
|---|---|---|---|
| $\hat{\beta}_1$ | Mean (SD) | 0.297 (0.284) | 0.291 (0.139) |
|  | Coverage |  |  |
|  | for $m = 10, 40$ | 0.960, 0.970 | 0.960, 0.970 |
| $\hat{\beta}_2$ | Mean (SD) | 0.247 (0.168) | 0.246 (0.083) |
|  | Coverage |  |  |
|  | for $m = 10, 40$ | 0.970, 0.990 | 0.970, 0.955 |
| Joint | Coverage |  |  |
|  | for $m = 10, 40$ | 0.975, 0.990 | 0.960, 0.955 |

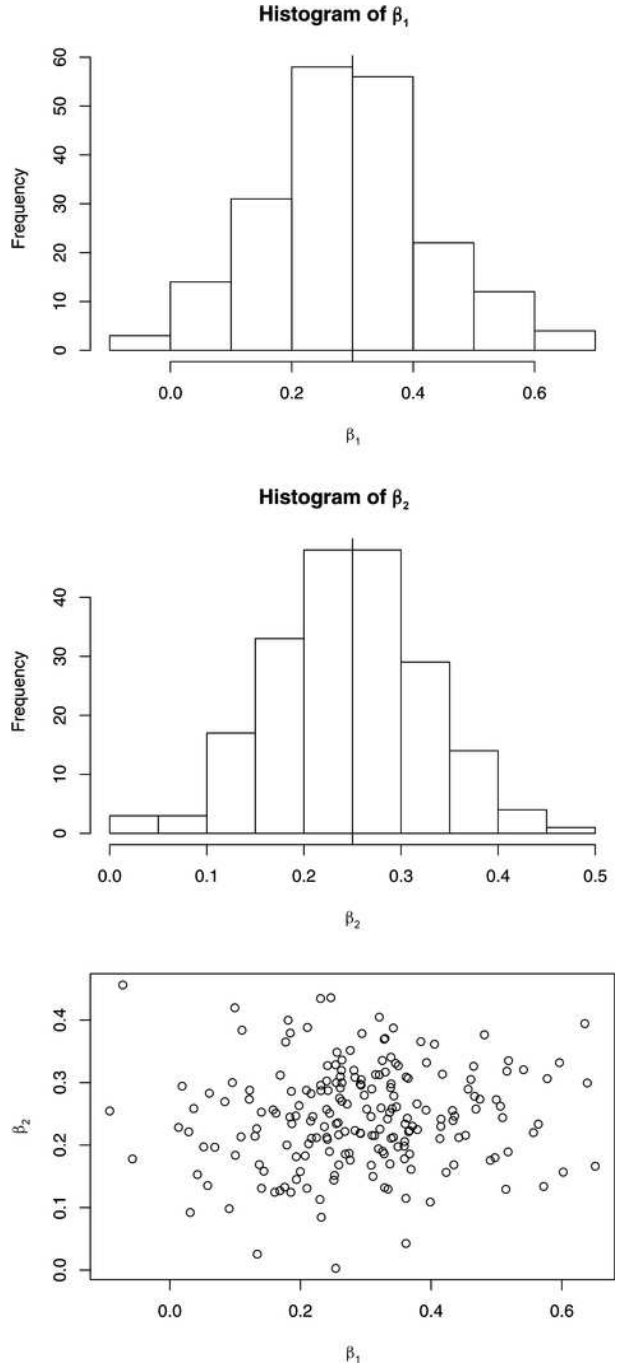

FIG. 1.  *Histogram of estimations of $\beta_1$ and $\beta_2$. Scatter plot of $\hat{\beta}_{1,n}$ versus $\hat{\beta}_{2,n}$. The sample size is 1600. Based on 200 replicates.*

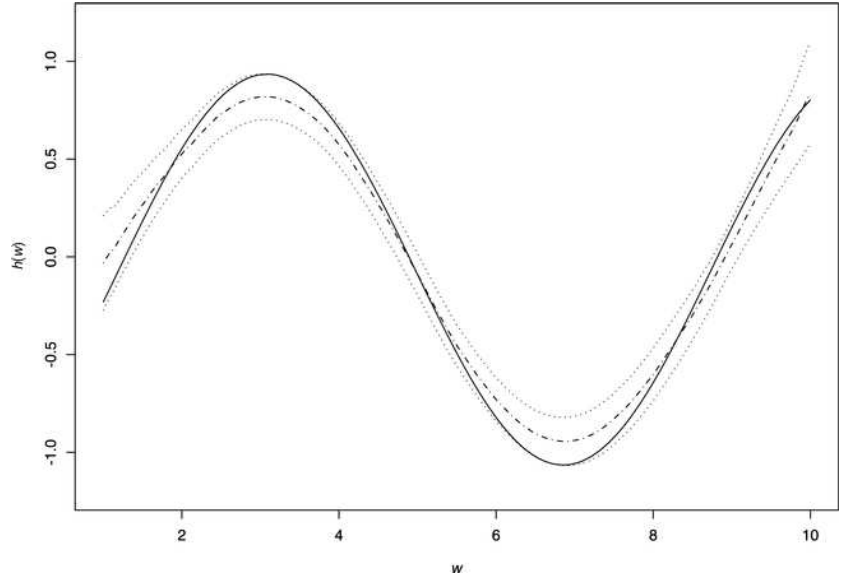

FIG. 2. *Estimate and pointwise confidence interval for h. The solid line is the true value. The dashed line is the estimated mean value. The dotted lines are the pointwise 95% confidence intervals. The sample size is 1600, based on 200 replicates.*

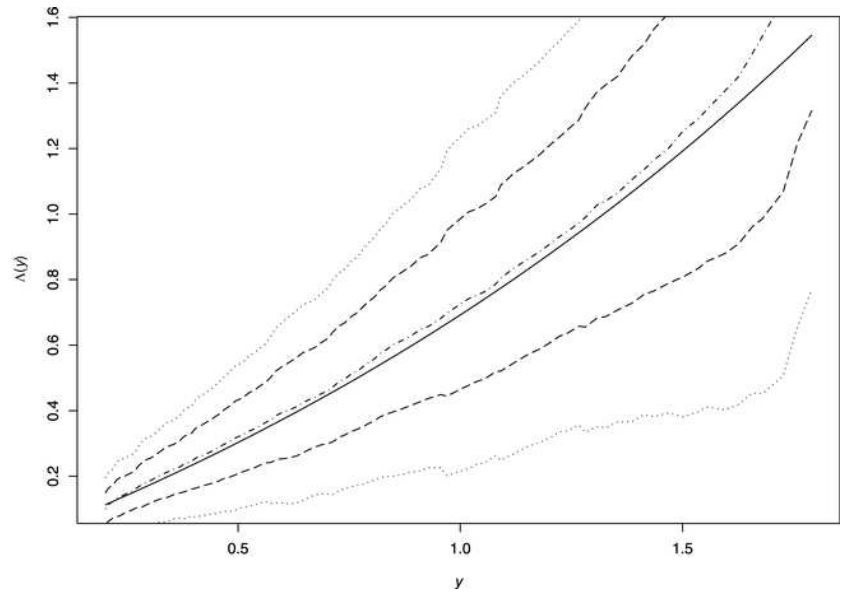

FIG. 3. *Estimation and pointwise confidence intervals for Λ. The solid line is the true value. The dot-dashed line by the solid line is the estimated mean value. The dashed lines are the mean value plus (minus) one standard deviation. The dotted lines are pointwise 95% confidence intervals. The sample size is 1600, based on 200 replicates.*

PROOF OF LEMMA 1.   By the inclusions contained in examples 3 and 4, and by Remark 2, it suffices to check example 1, example 3 for $\gamma \in (0, \infty)$, and example 5 for $\gamma \in (1, \infty)$. It is straightforward to verify B5(b) for each of these distributions. We now verify B5(d). For $F(u) = 1 - \exp[-e^u]$, $f(u) = e^u \exp[-e^u]$ and $\dot{f}(u) = -e^u(e^u - 1)\exp[-e^u]$. Thus, $f^2(u) - \dot{f}(u)F(u) = e^u \exp[-e^u](\exp[-e^u] + e^u - 1) > 0$ for all $u \in \mathbb{R}$, since $e^{-v} + v$ is strictly increasing on $(0, \infty)$. Also, $f^2(u) + \dot{f}(u)(1 - F(u)) = e^u \exp[-e^u] > 0$ for all $u \in \mathbb{R}$. Hence, B5(d) is satisfied for example 1. Similar arguments establish the condition for example 3. Consider now $F(u) = r_\gamma \int_{-\infty}^u e^{-|s|^\gamma} ds$ for $\gamma \in (1, \infty)$, where $r_\gamma \equiv \gamma[2\Gamma(1/\gamma)]^{-1}$. Since $f(u) = r_\gamma e^{-|u|^\gamma}$, $\dot{f}(u) = -\mathrm{sign}(u)\gamma r_\gamma |u|^{\gamma-1} e^{-|u|^\gamma}$ and, thus,

$$
\begin{aligned}
(10.1) \quad & f^2(u) + \dot{f}(u)(1 - F(u)) \\
& = r_\gamma^2 e^{-2|u|^\gamma} \left(1 - \frac{\mathrm{sign}(u)\gamma |u|^{\gamma-1} \int_{-\infty}^u e^{-|s|^\gamma} ds}{e^{-|u|^\gamma}}\right),
\end{aligned}
$$

which is clearly $> 0$ for all $u \in (-\infty, 0]$. Since $\int_v^\infty e^{-s^\gamma} ds < \gamma^{-1} v^{-(\gamma-1)} e^{-v^\gamma}$ for all $v \in (0, \infty)$, (10.1) $> 0$ for all $u \in (0, \infty)$. Similar techniques verify that $f^2(u) - \dot{f}(u)F(u) > 0$ for all $u \in \mathbb{R}$ and, thus, B5(d) is satisfied yet again.

Establishing condition B5(c) is more challenging. For $F(v) = 1 - \exp[-e^v]$, $F(F^{-1}(u) + s) = u^{e^s}$. Let $\eta_1 = \inf_{s \in K} e^s$, $\eta_2 = \sup_{s \in K} e^s$ and $\xi(s) = s^{\eta_1}$. Note that $\xi : [0, 1] \mapsto [0, 1]$ is increasing and isomorphic. Furthermore, $|u^{e^s} - v^{e^s}| = |\{\xi^{-1}\xi(u)\}^{e^s} - \{\xi^{-1}\xi(v)\}^{e^s}| \leq (\eta_2/\eta_1)|\xi(u) - \xi(v)|$, since $\inf_{s \in K} \eta_1^{-1} e^s = 1$ and $\sup_{s \in K} \eta_1^{-1} e^s = \eta_2/\eta_1$. Thus, (c) is satisfied for example 1 when $\alpha = 1$. For example 3 with $\gamma \in (0, \infty)$, $F(u) = 1 - [1 + \gamma e^u]^{-1/\gamma}$ and, thus, $F(F^{-1}(u) + s) = 1 - [1 - e^s + e^s(1-u)^{-\gamma}]^{-1/\gamma}$. Hence, after some derivation,

$$
\begin{aligned}
\frac{\partial}{\partial s} F(F^{-1}(u) + s) &= \frac{e^s}{\gamma}\{1 + e^s[(1-u)^{-\gamma} - 1]\}^{-1/\gamma-1}[(1-u)^{-\gamma} - 1] \\
&\leq \frac{e^s}{\gamma}\{1 + e^s[(1-u)^{-\gamma} - 1]\}^{-1/\gamma} \\
&\quad \times [1 - (1-u)^\gamma]\{(1-u)^\gamma + [1 - (1-u)^\gamma]e^s\}^{-1} \\
&\leq \frac{e^s \vee e^{-s}}{\gamma} \equiv c_*,
\end{aligned}
$$

which is uniformly bounded over $s \in K$. Thus, B5(c) is also satisfied in this instance, with $c = c_*$, $\alpha = 1$ and $\xi$ equal to the identity.

It appears to be quite difficult to establish B5(c) directly for example 5 with $\gamma \in (1, \infty)$, so we will use Lemma 2. In this case, we have that, for $u, \tau \geq 0$, $f(-u)\dot{f}(-u + \tau) - \dot{f}(-u)f(-u + \tau) \leq 0$ if and only if $0 \geq \mathrm{sign}(u - \tau) \times$ $|u - \tau| - u = -\tau$, but this last inequality is always true. Similarly, $f(u)\dot{f}(u - \tau) - \dot{f}(u)f(u - \tau) \geq 0$ if and only if $0 \leq -\mathrm{sign}(u - \tau)|u - \tau| + u = \tau$, but, again, this last inequality is always true whenever $\tau \geq 0$. Hence, condition (i) of Lemma 2 holds with $c_1^{(\tau)} = 0$ for all $\tau \in [0, \infty)$. Establishing condition (ii) is more challenging. Condition (ii) is trivially true when $\tau = 0$. Fix $\tau \in (0, \infty)$ and $\eta \in (0, 1/3)$. Let $G \colon [0, \infty) \mapsto [0, 1]$ be strictly increasing with $F(u) \leq G(u)$ for all $u > 0$. Then $u \geq G^{-1}(F(u))$ for all $u > 0$ and, thus, $1 - F(F^{-1}(1 - \epsilon) - \tau) \leq 1 - F(G^{-1}(1 - \epsilon) - \tau)$. Note that, for $u > 0$,

$$r_\gamma \int_u^\infty e^{-s^\gamma}\,ds = r_\gamma \int_u^\infty s^{\gamma-1} e^{-(1+\eta)s^\gamma}[s^{-(\gamma-1)}e^{\eta s^\gamma}]\,ds$$

$$\geq c_1 \int_u^\infty \gamma(1+\eta)s^{\gamma-1}e^{-(1+\eta)s^\gamma}\,ds$$

$$= c_1 e^{-(1+\eta)u^\gamma} \equiv 1 - G(u),$$

where $c_1 \in (0, \infty)$ does not depend on $u$ and the inequality follows since $\inf_{s>0}[s^{-(\gamma-1)}e^{\eta s^\gamma}] > 0$. Solving for $G(u_*) = 1 - \epsilon$, we obtain $u_* = [\log(c_1/\epsilon)/(1+\eta)]^{1/\gamma}$. Since it is also true that $1 - F(u) \leq c_2 e^{-u^\gamma}$ when $u > 0$ for some $c_2 \in (0, \infty)$ which does not depend on $u$,

$$
\begin{aligned}
(10.2) \quad 1 - F(F^{-1}(1 - \epsilon) - \tau) &\leq 1 - F(G^{-1}(1 - \epsilon) - \tau) \\
&\leq c_2 \exp\left(-\left\{\left[\frac{\log(c_1/\epsilon)}{1+\eta}\right]^{1/\gamma} - \tau\right\}^\gamma\right) \\
&= c_2 \exp\left\{-(1-\eta)\left[\frac{\log(c_1/\epsilon)}{1+\eta}\right] - q(\epsilon)\right\},
\end{aligned}
$$

where

$$q(\epsilon) \equiv \left[\frac{\log(c_1/\epsilon)}{1+\eta}\right]^{1/\gamma} - (1-\eta)\left[\frac{\log(c_1/\epsilon)}{1+\eta}\right]$$

is bounded below for all $\epsilon$ small enough. Hence, $1 - F(F^{-1}(1 - \epsilon) - \tau) \leq k_* \epsilon^{(1-\eta)/(1+\eta)}$ for some $k_* \in (0, \infty)$ not depending on $\epsilon$, for all $\epsilon$ small enough. Thus, condition (ii) is satisfied, and the lemma yields that B5(c) is satisfied in this setting. □

PROOF OF LEMMA 2. Fix a compact $K \subset \mathbb{R}$, and set $\tau = \sup\{|s| : s \in K\}$. Choose $\epsilon_1 \in (0, 1/3)$ so that condition (ii) of the lemma is satisfied for all $\epsilon \leq \epsilon_1$ and $F^{-1}(1 - \epsilon_1) \wedge [-F^{-1}(\epsilon_1)] > \tau + c_1^{(\tau)}$. Note that, for all $u \in [1 - \epsilon_1, 1]$ and $\rho \in [0, 1 - u]$,

$$
\begin{aligned}
\frac{\partial}{\partial s}[F(F^{-1}(\xi_*^{(\tau)}\{u + \rho\}) + s) &- F(F^{-1}(\xi_*^{(\tau)}\{u\}) + s))]|_{s=t} \\
&= f(F^{-1}(\xi_*^{(\tau)}\{u + \rho\}) + t) - f(F^{-1}(\xi_*^{(\tau)}\{u\}) + t) \leq 0,
\end{aligned}
$$

for all $t \in K$, which implies $F(F^{-1}(\xi_*^{(\tau)}\{u + \rho\}) + t) - F(F^{-1}(\xi_*^{(\tau)}\{u\}) + t) \le F(F^{-1}(\xi_*^{(\tau)}\{u + \rho\}) - \tau) - F(F^{-1}(\xi_*^{(\tau)}\{u\}) - \tau)$ for all $t \in K$. Arguing in a similar manner, we obtain for all $u \in [0, \epsilon_1]$ and $\rho \in [0, u]$ that $F(F^{-1}(\xi_*^{(\tau)}\{u\}) + t) - F(F^{-1}(\xi_*^{(\tau)}\{u - \rho\}) + t) \le F(F^{-1}(\xi_*^{(\tau)}\{u\}) + \tau) - F(F^{-1}(\xi_*^{(\tau)}\{u - \rho\}) + \tau)$ for all $t \in K$.

By condition (i) of Lemma 2, we have for all $\epsilon \in [0, \epsilon_1]$ that

$$\sup_{t \in K, u_1, u_2 \in [0, \epsilon_1] : |u_1 - u_2| \le \epsilon} |F(F^{-1}(\xi_*^{(\tau)}\{u_1\}) + t) - F(F^{-1}(\xi_*^{(\tau)}\{u_2\}) + t)|$$

$$\le \sup_{0 \le u_2 - \epsilon \le u_1 \le u_2 \le \epsilon_1} F(F^{-1}(\xi_*^{(\tau)}\{u_2\}) + \tau) - F(F^{-1}(\xi_*^{(\tau)}\{u_1\}) + \tau)$$

$$\le F(F^{-1}(\xi_*^{\tau}\{\epsilon\}) + \tau)$$

and

$$\sup_{t \in K, u_1, u_2 \in [1 - \epsilon_1, 1] : |u_1 - u_2| \le \epsilon} |F(F^{-1}(\xi_*^{(\tau)}\{u_1\}) + t) - F(F^{-1}(\xi_*^{\tau}\{u_2\}) + t)|$$

$$\le \sup_{1 - \epsilon_1 \le u_1 \le u_2 \le u_1 + \epsilon \le 1} F(F^{-1}(\xi_*^{(\tau)}\{u_2\}) - \tau) - F(F^{-1}(\xi_*^{(\tau)}\{u_1\}) - \tau)$$

$$\le 1 - F(F^{-1}(\xi_*^{(\tau)}\{1 - \epsilon\}) - \tau).$$

Since both $F$ and $F^{-1}$ have bounded derivatives on compacts,

$$\sup_{t \in K, u_1, u_2 \in [\epsilon_1, 1 - \epsilon_1] : |u_1 - u_2| \le \epsilon} |F(F^{-1}(\xi_*^{\tau}\{u_1\}) + t) - F(F^{-1}(\xi_*^{(\tau)}\{u_2\}) + t)| \le \tilde{c}\epsilon$$

for some $\tilde{c} \in (0, \infty)$. Hence,

$$\sup_{s \in K} \sup_{u, v \in [0, 1] : |u - v| \le \epsilon} |F(F^{-1}(\xi_*^{(\tau)}\{u\}) + s) - F(F^{-1}(\xi_*^{(\tau)}\{v\}) + s)| \le c\epsilon^{\alpha},$$

where $c = (1/\epsilon_1) \vee c_1^{(\tau)} \vee \tilde{c}$ and $\alpha = \alpha_\tau$. Since the compact set $K$ was arbitrary, the desired result follows with $\xi$ chosen so that $\xi^{-1} = \xi_*^{(\tau)}$.  $\square$

PROOF OF LEMMA 3.  Since $l(x; \beta, h, H) \le 0$, $\lambda_n > 0$ forces $\hat{h}_n \in \Im_\nu$ since otherwise $l_n^p(\hat{\beta}_n, \hat{h}_n, \hat{H}_n) = -\infty$. Since $B_0$ is bounded, $\hat{\beta}_n$ is obviously bounded. Since $\hat{h}_n$ is also bounded by assumption, $\lim_{H \downarrow -\infty} F(\hat{\beta}'_n z_{(1)} + \hat{h}_n(w_{(1)}) + H) = 0$ and $\lim_{H \uparrow \infty} F(\hat{\beta}'_n z_{(n)} + \hat{h}_n(w_{(n)}) + H) = 1$ and, thus, $l_n^p(\hat{\beta}_n, \hat{h}_n, \hat{H}_n) = -\infty$ if either $\hat{H}_n(v_{(1)}) = -\infty$ or $\hat{H}_n(v_{(n)}) = \infty$.  $\square$

PROOF OF THEOREM 1.  Define

$$l^*(x; \beta, h, H) \equiv \delta F(\beta' z + h(w) + H(v))$$

$$+ (1 - \delta)\{1 - F(\beta' z + h(w) + H(v))\}$$

and fix $\gamma \in (0,1)$. Note that since $l_n^p(\hat{\beta}_n, \hat{h}_n, \hat{H}_n) \geq l_n^p(\beta_0, h_0, H_0)$, $\lambda_n^2 J^2(\hat{h}_n) = O_p(1)$ and, thus, $J(\hat{h}_n) = O_p(n^{1/3})$. Also, $l_n(\hat{\beta}_n, \hat{h}_n, \hat{H}_n) \geq l_n(\beta_0, h_0, H_0) + O_p(n^{-2/3})$, which implies by the concavity of $s \mapsto \log(s)$ that

$$
\begin{aligned}
(10.3) \quad O_p(n^{-2/3}) &\leq \mathbb{P}_n \log\left[1 + \gamma \left\{ \frac{l^*(X; \hat{\beta}_n, \hat{h}_n, \hat{H}_n)}{l^*(X; \beta_0, h_0, H_0)} \right\} \right] \\
&\equiv \mathbb{P}_n \zeta(X; \hat{\beta}_n, \hat{h}_n, \hat{H}_n).
\end{aligned}
$$

These facts, combined with the result from Lemma 8 below that

$$
(10.4) \quad \mathcal{F}_0 = \{[1 + J(h)]^{-1} \zeta(X; \beta, h, H) : \beta \in \bar{B}_0, H \in \mathcal{M}_\mathbb{R}, h \in \Im_\nu\}
$$

is $P$-Donsker, where $\mathcal{M}_A$ is the collection of all nondecreasing functions mapping from $\mathbb{R}$ to $A \subset \mathbb{R}$, imply that $(\mathbb{P}_n - P)\zeta(X; \hat{\beta}_n' Z, \hat{h}_n(W), \hat{H}_n(V)) = O_p(n^{-1/6})$. Thus, $P\zeta(X; \hat{\beta}_n, \hat{h}_n, \hat{H}_n) \geq O_p(n^{-1/6})$. However, by the concavity of $s \mapsto \log(s)$, $P\zeta(X; \hat{\beta}_n, \hat{h}_n, \hat{H}_n) \leq 0$ and, thus, $P\zeta(X; \hat{\beta}_n, \hat{h}_n, \hat{H}_n) = O_p(n^{-1/6})$. Since $U_n(x) \equiv l^*(x; \hat{\beta}_n, \hat{h}_n, \hat{H}_n)/l^*(x; \beta_0, h_0, H_0)$ satisfies $0 \leq U_n(x) \leq m < \infty$ for some $m$ not depending on $n$ or $x$, Prohorov's theorem now implies for every subsequence $n'$ that there exists a further subsequence $n''$ so that $U_{n''}(X)$ converges in distribution to some $U(X)$ satisfying $P \log\{1 + \gamma(U(X) - 1)\} = 0$ and $PU(X) = 1$. But this implies $U(X) = 1$, almost surely, by the strict concavity of $s \mapsto \log\{1 + \gamma(s - 1)\}$. Since this result is true for every subsequence, we have that $P|U_n(X) - 1| = o_p(1)$. This now implies that

$$
\begin{aligned}
(10.5) \quad P\{F(\hat{H}_n(V) &+ \hat{\beta}_n' Z + \hat{h}_n(W)) \\
&- F(H_0(V) + \beta_0' Z + h_0(W))\}^2 = o_p(1).
\end{aligned}
$$

Expression (10.5) implies that

$$
P[\{(\hat{\beta}_n - \beta_0)'(Z - E[Z|V, W]) + c_n(V, W)\}^2 | V, W] = o_p(1),
$$

for almost surely all $V$ and $W$, where $c_n(V, W) \equiv (\hat{\beta}_n - \beta_0)' E[Z|V, W] + \hat{H}_n(V) - H_0(V) + \hat{h}_n(W) - h_0(W)$ is a sequence uncorrelated with $Z - E[Z|V, W]$. Condition A2 now implies that $\hat{\beta}_n - \beta_0 = o_p(1)$, and, furthermore, that

$$
(10.6) \quad P\{F(\hat{H}_n(V) + \hat{h}_n(W)) - F(H_0(V) + h_0(W))\}^2 = o_p(1).
$$

Let $\mathcal{V}$ be the set of all $V$ such that the distribution of $W$ given $V$ dominates the unconditional distribution of $W$. Condition A3, combined with (10.6), implies that for some $v \in \mathcal{V}$,

$$
\begin{aligned}
(10.7) \quad o_p(1) &= P[\{\hat{H}_n(V) + \hat{h}_n(W) - H_0(V) - h_0(W)\}^2 | V = v] \\
&\geq P_{W_1}\{\hat{h}_n(W_1) - h_0(W_1) - P_{W_2}\hat{h}_n(W_2)\}^2,
\end{aligned}
$$

where $P_{W_j}$ is the marginal probability measure of $W$ applied to $W_j$, $j = 1, 2$. Since $\{h/[1 + J(h)] : h \in \Im_\nu\}$ is Donsker (see Theorem 2.4 of [43]), $J(\hat{h}_n) = O_p(n^{1/3})$ and $P_{W_1} h_0(W_1) = 0$, we have that $P_{W_1} \hat{h}_n(W_1) = O_p(n^{-1/6})$. Thus, the last term in (10.7) implies

$$P_{W_1}\{\hat{h}_n(W_1) - h_0(W_1)\}^2 = o_p(1)$$

and, thus, by condition A1(b), $\|\hat{h}_n - h_0\|_2 = o_p(1)$. This now implies that $\|\hat{H}_n - H_0\|_{F,2} = o_p(1)$.  $\square$

PROOF OF THEOREM 2. Assume without loss of generality that $\Delta_{(1)} = 1$ and $\Delta_{(n)} = 0$ as discussed in Remark 6. Divide the observations into contiguous disjoint segments $M_k \subset \{1, \ldots, n\}$, $k = 1, \ldots, K$, where $1 = \min(M_1) < \max(M_1) = \min(M_2) - 1 < \min(M_2) < \cdots < \max(M_{k-1}) = \min(M_k) - 1 < \min(M_k) < \max(M_k) = n$, so that $\{\Delta_{(j)}, j \in M_k\}$ consists of all 1's followed by all 0's. Hence, there are at least two observations in each $M_k$, $k = 1, \ldots, K$. Note that $\hat{H}_n(V_{(i)}) = \hat{H}_n(V_{(j)})$ for all $i, j \in M_k$, $k = 1, \ldots, K$. To see this, suppose that it is not true, and let $j'$ be the index in $M_k$ which corresponds to the first time $\delta_{(j)} = 0$ over $j \in M_k$. Now the profile log-likelihood $l_n^p(\beta, h, \hat{H}_n)$ can be increased by replacing $\hat{H}_n$ with $\hat{H}_n^*$, where $\hat{H}_n^*(V_{(i)}) = \hat{H}_n(V_{(j')})$ for all $i < j', i \in M_k$ [since this would increase the value of $\log\{F(\hat{H}_n(V_{(i)}) + \hat{\beta}_n' Z_{(i)} + \hat{h}_n(W_{(i)}))\}$]. The profile log-likelihood will also be increased by setting $\hat{H}_n^*(V_{(i)}) = \hat{H}_n(V_{(j')})$ for all $i \geq j', i \in M_k$ [since this would lower $\log\{F(\hat{H}_n(V_{(i)}) + \hat{\beta}_n' Z_{(i)} + \hat{h}_n(W_{(i)}))\}$ and, hence, increase $\log\{1 - F(\hat{H}_n(V_{(i)}) + \hat{\beta}_n' Z_{(i)} + \hat{h}_n(W_{(i)}))\}$]. Thus, $l_n^p(\hat{\beta}_n, \hat{h}_n, \hat{H}_n^*) > l_n^p(\hat{\beta}_n, \hat{h}_n, \hat{H}_n)$, and the MLE $\hat{H}_n(V_{(i)})$ is therefore constant over the indices $i \in M_k$, $k = 1, \ldots, K$.

Define $v_0 \equiv H_0(l_v)$ and $p_0 \equiv F(v_0 - 2m) \wedge [1 - F(v_0)]$, where $m$ is the maximum possible value of $|\hat{\beta}_n' z_{(i)} + \hat{h}_n(w_{(i)})|$ over $1 \leq i \leq n$, and let $q(\delta, t) \equiv \delta \log(F(t)) + (1 - \delta) \log(1 - F(t))$. Note that condition 5(d) implies that $q(\delta, t)$ is strictly convex over $t \in \mathbb{R}$ for $\delta \in \{0, 1\}$. Accordingly, $\hat{H}_n(V_{(1)})$ is piecewise constant for all indices in $M_k$, $k = 1, \ldots, k^*$, for some $k^* \leq K$ and, thus, for any $\epsilon \in (0, p_0)$,

$$P\{\hat{H}_n(V_{(1)}) \leq F^{-1}(\epsilon)\}$$

$$\leq P\left\{\inf_{1 \leq j \leq K}\left[\operatorname*{arg\,max}_{a \in \mathbb{R}}\left(\sum_{l=1}^{j}\sum_{i \in M_l} \Delta_{(i)} \log F(a + \hat{\beta}_n' Z_{(i)} + \hat{h}_n(W_{(i)}))\right.\right.\right.$$

$$\left.\left.\left. + (1 - \Delta_{(i)}) \log\{1 - F(a + \hat{\beta}_n' Z_{(i)} + \hat{h}_n(W_{(i)}))\}\right)\right]\right.$$

$$(10.8) \hspace{4cm} \leq F^{-1}(\epsilon) \Big\}$$

$$\leq P\Big\{ \sup_{a \leq F^{-1}(\epsilon), 1 \leq j \leq K} \Big[ \sum_{l=1}^{j} \sum_{i \in M_l} \Big( \Delta_{(i)} \log\Big\{ \frac{F(a + \hat{\beta}'_n Z_{(i)} + \hat{h}_n(W_{(i)}))}{F(v_0 - m + \hat{\beta}'_n Z_{(i)} + \hat{h}_n(W_{(i)}))} \Big\}$$

$$+ (1 - \Delta_{(i)}) \log\Big\{ \frac{1 - F(a + \hat{\beta}'_n Z_{(i)} + \hat{h}_n(W_{(i)}))}{1 - F(v_0 - m + \hat{\beta}'_n Z_{(i)} + \hat{h}_n(W_{(i)}))} \Big\} \Big) \Big] \geq 0 \Big\}$$

$$\leq P\Big\{ \sup_{a \leq F^{-1}(\epsilon)} \sup_{1 \leq j \leq n} \Big[ \sum_{l=1}^{j} \sum_{i \in M_l} (\Delta_{(i)} \log F(a+m) - \log p_0) \Big] \geq 0 \Big\}$$

$$\leq P\Big\{ \inf_{1 \leq j \leq n} \bar{Q}_{(j)} \leq \frac{\log(1/p_0)}{\log[1/F(F^{-1}(\epsilon) + m)]} \Big\},$$

where $\bar{Q}_{(j)} = j^{-1} \sum_{i=1}^{j} \Delta_{(i)}$.

Let $\bar{Q}^*_j = (j+1)^{-1}[1 + \sum_{i=1}^{j-1} Q^*_i]$, where $Q^*_1, Q^*_2, \dots$ are i.i.d. Bernoulli with probability of success $F(H_0(l_v) - m)$. Then (10.8) is bounded above by

$$(10.9) \hspace{2cm} P\Big\{ \inf_{j \geq 1} \bar{Q}^*_j \leq \frac{\log(1/p_0)}{\log[1/F(F^{-1}(\epsilon) + m)]} \Big\}.$$

For every $\tau \in (0, F(H_0(l_v) - m))$, the strong law of large numbers yields that $N_\tau = \sup\{j : j^{-1} \sum_{i=1}^{j} Q^*_i \leq \tau\}$ is a bounded random variable. This now implies that the probability in (10.8) can be made arbitrarily small by taking $\epsilon$ small enough. Hence, $|\hat{H}_n(V_{(1)})| = O_p(1)$. The proof that $|\hat{H}_n(V_{(n)})| = O_p(1)$ is obtained in a virtually identical manner after reversing the order of the indices. $\square$

PROOF OF THEOREM 3. By the isotonic regression results in Section II.1.1 of [19],

$$\hat{H}_n(V_{(1)}) = \min_{1 \leq k \leq n} \frac{\sum_{i=1}^{k} \Delta_{(i)}}{k} \leq \frac{1}{M_n},$$

where $M_n = \max\{j \leq n : \sum_{i=1}^{j} \Delta_{(i)} = 1\}$. In this setting $\Delta_{(1)} = 1$ and $\Delta_{(n)} = 0$ almost surely (by assumption); but $\Delta_{(2)}, \dots, \Delta_{(n-1)}$ conditional on $V_{(2)}, \dots, V_{(n-1)}$ are independent Bernoullis with probabilities of success $G_0(V_{(i)})$, $i = 2, \dots, n-1$. Thus, $M_n$ is bounded below in probability by $M^*_n \equiv 1 + \max\{j \leq n-1 : \sum_{i=1}^{j} \Delta^*_i = 0\}$, where the $\Delta^*_i$'s are i.i.d. Bernoullis with probability of success $G_0(u_v) < 1$. Let $M^* \equiv \lim_{n \to \infty} M^*_n$. Since for any $k < \infty$, $P\{M^* \geq k\} = \{1 - G_0(u_v)\}^{k-1} > 0$, we now have for any $\eta > 0$ that $\liminf_{n\to\infty} P\{\widehat{G}_n(V_{(1)}) \le \eta\} \ge P\{1/M^* \le \eta\} > 0$. Thus, (i) follows. The same argument can be used to verify (ii) after noting that

$$1 - \widehat{H}_n(V_{(n)}) = 1 - \max_{1 \le k \le n} \frac{\sum_{i=k}^n \Delta_{(i)}}{n-k+1} = \min_{1 \le k \le n} \frac{\sum_{i=1}^k \{1 - \Delta_{(n-i+1)}\}}{k}. \quad \square$$

PROOF OF THEOREM 4. We make use of the following technical tools.

T1. Denote $\Theta = \{\theta : \theta = g + H, g \in \mathfrak{G}, H \in \Xi\}$, where $\Xi = \{H : H$ is a nondecreasing function and $-\infty < M_1 \le H \le M_2 < \infty\}$, for constants $M_1$ and $M_2$, and where $\mathfrak{G} = \{g : g = \beta'z + h(w) : \beta \in \bar{B}_0, |h| \le c_0, J(h) < \infty\}$. The arguments in Lemma 8 yield that $\log N_{[\cdot]}(\epsilon, \Theta/(1 + J(h)), P) \le A_1\epsilon^{-1}$, where $A_1$ is a constant.

T2. (Theorem in [43], page 79.) Consider a uniformly bounded class of functions $\mathcal{G}$, with $\sup_{g \in \mathcal{G}} |g - g_0|_\infty < \infty$ and $\log N_{[\cdot]}(\epsilon, \mathcal{G}, P) \le A\epsilon^{-\alpha}$ for all $\epsilon > 0$, and where $\alpha \in (0, 2)$.

Then for $\delta_n = n^{-1/(2+\alpha)}$,

$$\sup_{g \in \mathcal{G}} \frac{|(\mathbb{P}_n - P)(g - g_0)|}{\|g - g_0\|_2^{1-\alpha/2} \vee \sqrt{n}\delta_n^2} = O_p(n^{-1/2}),$$

where $\|\cdot\|_2$ is the $L_2(P)$ norm.

Denote $\tilde{\theta}_n(x) \equiv \hat{\beta}_n'z + \hat{h}_n(w) + \widetilde{H}_n(v)$, $\theta_0(x) \equiv \beta_0'z + h_0(w) + H_0(v)$ and $q(\delta, t) \equiv \delta \log(F(t)) + (1 - \delta)\log(1 - F(t))$. Then $l(\beta, h, H)(x) = q(\delta, \theta(x))$. Denote the second-order derivative of $q$ as $-m(\delta, t) \equiv (\partial^2/\partial t^2)q(\delta, t)$.

Since $(g_0, H_0)$ maximizes the expectation of the log-likelihood function, we have

$$(10.10) \qquad P[l(g_0, H_0) - l(\hat{g}_n, \widetilde{H}_n)] = P\left[\frac{m(\delta, t^*)}{2}(\tilde{\theta} - \theta_0)^2\right],$$

where $t^*(X)$ is on the line segment between $\tilde{\theta}_n(X)$ and $\theta_0(X)$. From the compactness of $\tilde{\theta}_n$ and $\theta_0$, as given in the assumptions and as a consequence of Theorem 2,

$$(10.11) \qquad \exists \varepsilon_1, \varepsilon_2 : 0 < \varepsilon_1 < \varepsilon_2 < \infty \quad \text{and} \quad \varepsilon_1 < m(\delta, t^*) < \varepsilon_2 \qquad \text{a.s.}$$

Combining (10.10) and (10.11),

$$(10.12) \qquad \varepsilon_1\|\tilde{\theta}_n - \theta_0\|_2^2 \le P[l(g_0, H_0) - l(\hat{g}_n, \widetilde{H}_n)] \le \varepsilon_2\|\tilde{\theta}_n - \theta_0\|_2^2.$$

The penalized MLE estimators satisfy

$$(10.13) \qquad \begin{aligned} \lambda_n^2 J^2(\hat{h}_n) &\le \lambda_n^2 J^2(h_0) + \mathbb{P}_n[l(\hat{\theta}_n) - l(\tilde{\theta}_n)] \\ &\quad + (\mathbb{P}_n - P)[l(\tilde{\theta}_n) - l(\theta_0)] + P[l(\tilde{\theta}_n) - l(\theta_0)], \end{aligned}$$

where $\hat{\theta}_n(x) \equiv \hat{\beta}'_n z + \hat{h}_n(w) + \widehat{H}_n(v)$. However, if we let $k$ and $l$ be as defined in Remark 7, then

$$\mathbb{P}_n[l(\hat{\theta}_n) - l(\tilde{\theta}_n)] = n^{-1}\left[\sum_{i=1}^{k-1} \log\{1 - F(\hat{g}_n(x_{(i)}) + \widehat{H}_n(y_{(k)}))\}\right.$$
$$\left. + \sum_{i=l+1}^{n} \log F(\hat{g}_n(x_{(i)}) + \widehat{H}_n(y_{(l)}))\right]$$
$$= O_p(n^{-1})$$

by arguments given in the proof of Theorem 2. This, combined with (10.12) and (10.13), yields

$$(10.14) \quad \begin{aligned} \lambda_n^2 J^2(\hat{h}_n) + \varepsilon_1\|\tilde{\theta}_n - \theta_0\|_2^2 \\ \leq \lambda_n^2 J^2(h_0) + (\mathbb{P}_n - P)[l(\tilde{\theta}_n) - l(\theta_0)] + O_p(n^{-1}). \end{aligned}$$

Combined with the results in T1 and T2 (for $\alpha = 1$), (10.14) gives us

$$\begin{aligned} \lambda_n^2 J^2(\hat{h}_n) + \varepsilon_1\|\tilde{\theta}_n - \theta_0\|_2^2 \\ \leq \lambda_n^2 J^2(h_0) + O_p(n^{-1/2})(1 + J(\hat{h}_n))(\|\tilde{\theta}_n - \theta_0\|_2^{1/2} \vee n^{-1/6}). \end{aligned}$$

Thus, we conclude

$$(10.15) \quad \begin{aligned} \lambda_n^2 J^2(\hat{h}_n) \leq \lambda_n^2 J^2(h_0) \\ + O_p(n^{-1/2})(1 + J(\hat{h}_n))(\|\tilde{\theta}_n - \theta_0\|_2^{1/2} \vee n^{-1/6}), \end{aligned}$$

as well as

$$(10.16) \quad \begin{aligned} \varepsilon_1\|\tilde{\theta}_n - \theta_0\|_2^2 \leq \lambda_n^2 J^2(h_0) \\ + O_p(n^{-1/2})(1 + J(\hat{h}_n))(\|\tilde{\theta}_n - \theta_0\|_2^{1/2} \vee n^{-1/6}). \end{aligned}$$

A few further calculations give us $J(\hat{h}_n) = O_p(1)$ and $\|\tilde{\theta}_n - \theta_0\|_2 = O_p(n^{-1/3})$. Following arguments similar to those used in the proof of Theorem 1, we conclude $\|\hat{\beta}_n - \beta_0\| + d((\hat{h}_n, \tilde{H}_n), (h_0, H_0)) = O_p(n^{-1/3})$. Moreover, since $J(\hat{h}_n) = O_p(1)$, we can now conclude uniform consistency of $\hat{h}_n$ as discussed in Remark 8. □

PROOF OF LEMMA 4. Note that the model assumptions ensure that $D_0$ is bounded above and below on $[l_v, u_v] \times [a, b]$. Let $\mathcal{S}_1$ be the the class of functions $g : [l_v, u_v] \mapsto \mathbb{R}$ with $E[g^2(V)D_0(V, W)] < \infty$, and let $\mathcal{S}_2$ be the class of functions $g : [a, b] \mapsto \mathbb{R}$ with $E[g^2(W)D_0(V, W)] < \infty$. Because $D_0$ is bounded above and below, it is not hard to show that the score spaces $\mathcal{S}_0^* \equiv \{g(V, W)Q_{\psi_0}(X) : g \in \mathcal{S}_0\}$ and $\mathcal{S}_1^* \equiv \{g(V)Q_{\psi_0}(X) : g \in \mathcal{S}_1\}$ are closed in $L_2(P)$, and that the $L_2(P)$ closure of the score space $\{g(W)Q_{\psi_0}(X) : g \in \Im_{\nu,0}\}$ is $\mathcal{S}_2^* \equiv \{g(W)Q_{\psi_0}(X) : g \in \mathcal{S}_2\}$. The reason we can drop the requirement $Eg(W) = 0$ in the latter case is that $E[Q_{\psi_0}(X)|V,W] = 0$. We also note that now both $\mathcal{S}_1^*$ and $\mathcal{S}_2^*$ are closed subspaces of $\mathcal{S}_0^*$. We can also see that, for any $g(V,W)Q_{\psi_0}(X) \in \mathcal{S}_0^*$, $(\Pi_1 g)(V)Q_{\psi_0}(X)$ is the projection onto $\mathcal{S}_1^*$ and $(\Pi_2 g)(W)Q_{\psi_0}(X)$ is the projection onto $\mathcal{S}_2^*$.

Define a new score space $\mathcal{S}_3^* \equiv \{[g(V) + h(W)]Q_{\psi_0}(X) : g \in \mathcal{S}_1, h \in \mathcal{S}_2\}$. Since $D_0$ is bounded below and $V$ and $W$ are independent, there exists a constant $c > 0$ such that for all $g \in \mathcal{S}_1$ with $E[g(V)] = 0$ and all $h \in \mathcal{S}_2$,

$$E\{[g(V) + h(W)]^2 Q_{\psi_0}^2(X)\} \geq cE[g(V) + h(W)]^2 \geq cEg^2(V) + cEh^2(W).$$

Thus, $\mathcal{S}_3^*$ is also a closed subspace of $\mathcal{S}_0^*$. This means that we can use the alternating projections theorem (Theorem A.4.2 of [6]) to establish that there exist a $q' \in \mathcal{S}_1$ and an $h' \in \mathcal{S}_2$ so that $D^*(V,W) - q'(V) - h'(W) = \lim_{m\to\infty}[(1-\Pi_2)(1-\Pi_1)]^m D^*$ and

$$(10.17) \qquad E\{[D^*(V,W) - q'(V) - h'(W)][q(V) + h(W)]Q_{\psi_0}^2(X)\} = 0,$$

for all $q \in \mathcal{S}_1$ and $h \in \mathcal{S}_2$. If we can show that

$$(10.18) \qquad \int_a^b \left[\left(\frac{\partial}{\partial w}\right)^\nu h'(w)\right]^2 dw < \infty$$

and that $h^* = h'$, then the first expression in D1 will hold for all $h \in \Im_{\nu,0}$. This last assertion follows by setting $q(V) = -(\Pi_1 h)(V)$ in (10.17) and noting that $q'(V)Q_{\psi_0}(X)$ is uncorrelated with $[h(W) - (\Pi_1 h)(V)]Q_{\psi_0}(X)$ for any $h \in \mathcal{S}_2 \supset \Im_{\nu,0}$.

We first establish that $h' = h^*$. For each $k \geq 0$, let $q_k \in \mathcal{S}_1$ and $h_k \in \mathcal{S}_2$ be defined by the equation

$$(10.19) \qquad D^*(V,W) - q_k(V) - h_k(W) = [(1-\Pi_2)(1-\Pi_1)]^k D^*.$$

To see that this makes sense, begin with $q_0 = h_0 = 0$, and note that, for any $k \geq 0$, $(1-\Pi_1)[D^*(V,W) - q_k(V) - h_k(W)] = D^* - (\Pi_1 D^*)(V) + (\Pi_1 h_k)(V) - h_k(W)$ and $(1-\Pi_2)[D^* - q_{k+1}(V) - h_k(W)] = D^* - q_{k+1}(V) - (\Pi_2 D^*)(W) + (\Pi_2 q_{k+1})(W)$. Thus, by setting $q_{k+1}(V) = (\Pi_1 D^*)(V) - (\Pi_1 h_k)(V)$ and

$$\begin{aligned}
(10.20) \qquad h_{k+1}(W) &= (\Pi_2 D^*)(W) - (\Pi_2 q_{k+1})(W) \\
&= [\Pi_2(1-\Pi_1)D^*](W) + (\Pi_2 \Pi_1 h_k)(W),
\end{aligned}$$

we have a method of defining $q_k$ and $h_k$ which is consistent with (10.19). By solving the recursive formula in (10.20), we obtain that $h_k(W) = [\sum_{j=0}^{k-1}(\Pi_2\Pi_1)^j] \times \Pi_2(1-\Pi_1)D^*$ for any $k \geq 1$. Thus, $h'$ in (10.17) is the limiting value of $h_k$, as $k \to \infty$, where the limit is in $\mathcal{S}_2$ since $\mathcal{S}_3^*$ is closed. But this is precisely how $h^*$ is defined. Thus, $h' = h^*$, and the limit in the definition of $h^*$ is well defined.

Now we will establish (10.18). Recall from above that $q' \in L_2(V)$. Note also that the above recursive arguments imply that $h^* = \Pi_2 D^* - \Pi_2 q'$. Thus, if we let $P_1$ be the probability measure for $V$,

$$\begin{aligned}
\|h^*\|_{P,\nu} &= \|\Pi_2 D^* - \Pi_2 q'\|_{P,\nu} \\
&\leq \|\Pi_2 D^*\|_{P,\nu} + \left\| \int_{l_v}^{u_v} q'(v) S_2(v,w) \, dP_1(v) \right\|_{P,\nu} \\
&\leq (E[D_2^{*(\nu)}(W)]^2)^{1/2} + \left( \int_{l_v}^{u_v} [q'(v)]^2 \, dP_1(v) \times E[S_2^{(\nu)}(V,W)]^2 \right)^{1/2} \\
&< \infty
\end{aligned}$$

by the boundedness assumptions on $D_2^{*(\nu)}$ and $S_2^{(\nu)}$. Since the density of $W$ is bounded below, we now have that $\int_a^b [(\partial/(\partial w))^\nu h^*(w)]^2 \, dw < \infty$. Thus, the first part of D1 is established and all that remains is to establish the required differentiability of $\tilde{q}$.

Because of the assumptions about $\dot{D}_1^*$, all that remains is establishing that $(\Pi_1 h^*)(v)$ has a derivative which is uniformly bounded on $[l_v, u_v]$. Letting $p_2(w)$ be the density of $W$, we have that $(\Pi_1 h^*)(v) = \int_a^b \tilde{h}(w) R_1(v,w) p_2(w) \, dw$. Thus,

$$|(\Pi_1 h^*)(v)| \leq \left( \int_a^b \tilde{h}^2(w) p_2(w) \, dw \times E[\dot{R}_1(v,W)]^2 \right)^{1/2},$$

and the desired result follows from the assumptions on $\dot{R}_1$. This completes the proof. $\square$

PROOF OF THEOREM 5. The information calculation is based on the non-orthogonal projection approach discussed by Sasieni [37]. The log-likelihood function takes the form

$$\begin{aligned}
l(x; \beta, h, H) &= \delta \log\{F[\beta'z + h(w) + H(v)]\} \\
&\quad + (1 - \delta) \log\{1 - F[\beta'z + h(w) + H(v)]\}.
\end{aligned}$$

The score function for $\beta$ is simply the derivative of the log-likelihood with respect to $\beta$, which is

$$\dot{l}_\beta = \delta Z \frac{f(\theta_\psi)}{F(\theta_\psi)} + (1 - \delta) Z \frac{-f(\theta_\psi)}{1 - F(\theta_\psi)} = Z f(\theta_\psi) \left( \frac{\delta}{F(\theta_\psi)} - \frac{1 - \delta}{1 - F(\theta_\psi)} \right),$$

where $\theta_\psi(x) \equiv \beta'z + h(w) + H(v)$.

Assume $h_\eta(w) = h(w) + \eta\xi(w)$, where $\xi \in \Im_{\nu,0}$. Then $(\partial/\partial\eta) h_\eta(w) = \xi(w)$. Thus, the score operator for $h(w)$ is

$$\dot{l}_h(\xi)(x) = \xi(w) f(\theta_\psi(x)) \left( \frac{\delta}{F(\theta_\psi(x))} - \frac{1 - \delta}{1 - F(\theta(x))} \right).$$

For the $H$ part, assume $(\partial/\partial\eta)H_\eta(v) = a(v)$, where $a \in L_2(V)$, and where $L_2(V)$ is the set of functions of the random variable $V$ which are square-integrable. Then

$$\dot{l}_H(a)(x) = a(v)f(\theta_\psi(x))\left(\frac{\delta}{F(\theta_\psi(x))} - \frac{1-\delta}{1-F(\theta_\psi(x))}\right).$$

*Step* 1. As suggested by Sasieni [37], we first project $\dot{l}_\beta(X)$ onto the space generated by $\dot{l}_H(X)$. We will need to find a function $a^* \in L_2(V)$ so that $\dot{l}_\beta - \dot{l}_H(a^*) \perp \dot{l}_H(a)$ for all $a \in L_2(V)$, which is equivalent to requiring

$$(10.21) \qquad E[(Z - a^*(V))a(V)Q_\psi^2(X)] = 0$$

for all $a \in L_2(V)$, where $Q_\psi$ is as defined in assumption D1. Since $E[(Z - a^*(V))a(V)Q_\psi^2(X)] = E[a(V)E((Z - a^*(V))Q_\psi^2(V)|V)] = 0$, then if $E((Z - a^*)Q^2|V) = 0$ almost surely, (10.21) will definitely be true. Thus, we can conclude that

$$(10.22) \qquad a^*(v) = \frac{E(ZQ_\psi^2(X)|V = v)}{E(Q_\psi^2(X)|V = v)}.$$

Hence, we have

$$(10.23) \qquad \dot{l}_\beta(X) - \dot{l}_H(a^*)(X) = \left(Z - \frac{E(ZQ_\psi^2(X)|V)}{E(Q_\psi^2|V)}\right)Q_\psi(V).$$

*Step* 2. We next project $\dot{l}_h(\xi)$ onto the space generated by $\dot{l}_H$, using similar calculations, to obtain a $b_* \in L_2(V)$ so that

$$(10.24) \quad \dot{l}_h(X) - \dot{l}_H(b^*)(X) = Q_\psi(X)\left\{\xi(W) - \frac{E(\xi(W)Q_\psi^2(X)|V)}{E(Q_\psi^2(X)|V)}\right\}.$$

*Step* 3. Next, we project the space generated by $\dot{l}_\beta - \dot{l}_H(a^*)$ onto the space generated by $\dot{l}_h - \dot{l}_H(b^*)$, which is equivalent to finding $\tilde{h} \in \Im_{\nu,0}$ such that

$$E\left\{\left(\left[Z - \frac{E(ZQ_\psi^2(X)|V)}{E(Q_\psi^2(X)|V)}\right] - \left[\tilde{h}(W) - \frac{E(\tilde{h}(W)Q_\psi^2(X)|V)}{E(Q_\psi^2(X)|V)}\right]\right)\right.$$

$$\left.\times \left[h(W) - \frac{E(h(W)Q_\psi^2(X)|V)}{E(Q_\psi^2(X)|V)}\right]Q_\psi^2(X)\right\} = 0$$

for all $h \in \Im_{\nu,0}$. That this is equivalent to the first conditional expectation in condition D1 follows from the fact that $r(V)Q_{\psi_0}(X)$ is uncorrelated with

$$\left[h(W) - \frac{E(h(W)Q_{\psi_0}^2(X)|V)}{E(Q_{\psi_0}(X)|V)}\right]Q_{\psi_0}(X)$$

for any $r \in L_2(V)$ and any $h \in L_2(W)$.

This proves Theorem 5. □

PROOF OF THEOREM 6. The proof of asymptotic normality and efficiency is based on Theorem 7 below, which is a modification of van der Vaart's Theorem 25.54 in [47].

As shown in Theorem 5, the efficient score function for $\beta$ takes the form $[z - \tilde{h}(w) - \tilde{q}(v)]Q_{\psi_0}(x)$, where $\tilde{h}$, $\tilde{q}$ and $Q_\psi$ are as defined in assumption D1. Formally, this function is the derivative at $t = 0$ of the log-likelihood function evaluated at $(\beta_0 + t, h_0 - t\tilde{h}, H_0 - t\tilde{q})$. However, the last coordinate of the latter path may not define a nondecreasing function for every $t$ in a neighborhood of 0 and, hence, cannot be used to obtain a stationary equation for the maximum likelihood estimator. To overcome this difficulty, we will replace the efficient score with an approximation based on an approximately least-favorable submodel.

For $t \in \mathbb{R}^d$, define $H_t(v) \equiv H(v) - t'\tilde{q}[H_0^{-1}(H_0(a) \vee [H(v) \wedge H_0(b)])]$. Then for $t$ close enough to zero, $H_t$ defines a nondecreasing function, since $v \mapsto \tilde{q}(H_0^{-1}(v))$ is Lipschitz continuous on $[H_0(a), H_0(b)]$ as a consequence of condition D1 and the assumed differentiability of $H_0^{-1}$. Now plug $(\beta + t, h_0 - t\tilde{h}, H_t)$ into the log-likelihood function and differentiate with respect to $t$ at $t = 0$. We then get the score function $\tilde{k}_\psi(X) \equiv (Z - \tilde{h}(W) - \tilde{q}[H_0^{-1}(H_0(a) \vee [H(V) \wedge H_0(b)])])Q_\psi(X)$, for which $\tilde{k}_{\psi_0}$ is the efficient score for $\beta$ at $\psi_0$.

We now have the following results.

1. The model is differentiable in quadratic mean with respect to $\beta$ at $(\beta_0, h_0, H_0)$.

2. As shown in Theorem 5, the efficient information matrix is nonsingular.

3. Note that, by the uniform consistency of $\hat{h}_n$, we have for large enough $n$ that $\hat{h}_n$ is in the interior of $\Im_\nu^{c_0}$ with high probability and satisfies $\mathbb{P}_n \hat{h}_n(W) = 0$. Hence, the derivative in the direction $\tilde{h} - \mathbb{P}_n \tilde{h}(W)$ of the log-likelihood will be zero for large enough $n$. This implies that

$$\mathbb{P}_n \tilde{k}_{\hat{\beta}_n, \hat{h}_n, \hat{H}_n} - \mathbb{P}_n[\tilde{h}]\mathbb{P}_n[Q_{\hat{\beta}_n, \hat{h}_n, \hat{H}_n}] = 0.$$

   Since $|\mathbb{P}_n(Q_{\hat{\beta}_n, \hat{h}_n, \hat{H}_n} - Q_{\hat{\beta}_n, \hat{h}_n, \tilde{H}_n})| = O_p(n^{-1})$ by arguments given in the proof of Theorem 2 and since $\mathbb{P}_n Q_{\hat{\beta}_n, \hat{h}_n, \tilde{H}_n} = o_p(1)$, we now have that $\mathbb{P}_n \tilde{k}_{\hat{\beta}_n, \hat{h}_n, \hat{H}_n} = o_p(n^{-1/2})$. Hence, also $\mathbb{P}_n \tilde{k}_{\hat{\beta}_n, \hat{h}_n, \tilde{H}_n} = o_p(n^{-1/2})$.

4. $(\hat{\beta}_n, \hat{h}_n, \tilde{H}_n)$ is consistent for $(\beta_0, h_0, H_0)$ (in an $L_2$ sense for $\tilde{H}_n$) and asymptotically bounded.

5. Since $J(\hat{h}_n) = O_p(1)$, and by results given in the proof of Theorem 5 and the Lipschitz continuity of the function $\psi \mapsto \tilde{k}_\psi$, we can see that there exists a neighborhood of $(\beta_0, h_0, H_0)$ such that the functions $\tilde{k}_{\hat{\beta}_n, \hat{h}_n, \tilde{H}_n}$ belong to a $P_{\beta_0, h_0, H_0}$ Donsker class with a square integrable envelope function with high probability for large enough $n$.

6. Denote $\zeta \equiv (h, H)$, with $\zeta_0 \equiv (h_0, H_0)$ and $\tilde{\zeta}_n \equiv (\hat{h}_n, \widetilde{H}_n)$. Let $B_{\beta, \zeta}$ be the score operator for $\zeta$ under the assumed model, but without the requirement that $P_{\beta, \zeta} h(W) = 0$. From the orthogonality of $\tilde{k}_{\beta, \zeta}$ and $B_{\beta, \zeta} D$ for any $D \in \Im_\nu \times L_2(V)$, we can now write

$$
\begin{aligned}
P_{\hat{\beta}_n, \zeta_0} \tilde{k}_{\hat{\beta}_n, \tilde{\zeta}_n} &= (P_{\hat{\beta}_n, \zeta_0} - P_{\hat{\beta}_n, \tilde{\zeta}_n})(\tilde{k}_{\hat{\beta}_n, \tilde{\zeta}_n} - \tilde{k}_{\beta_0, \zeta_0}) \\
(10.25) \qquad &\quad - \int \tilde{k}_{\beta_0, \zeta_0}[p_{\beta_0, \tilde{\zeta}_n} - p_{\beta_0, \zeta_0} - B_{\beta_0, \zeta_0}(\tilde{\zeta}_n - \zeta_0)p_{\beta_0, \zeta_0}]\, d\mu \\
&\quad + \int \tilde{k}_{\beta_0, \zeta_0}[p_{\beta_0, \tilde{\zeta}_n} - p_{\beta_0, \zeta_0} - p_{\hat{\beta}_n, \tilde{\zeta}_n} + p_{\hat{\beta}_n, \zeta_0}]\, d\mu,
\end{aligned}
$$

where $\mu$ is a suitable dominating measure. To verify the "no bias" condition (10.27) of Theorem 7 below, we use the decomposition (10.25). By the boundedness of the second derivative of $\log p_{\beta, \zeta}$ in a neighborhood of $(\beta_0, \zeta_0)$, the first term on the right-hand side of (10.25) is bounded by $O_p(1)d(\tilde{\zeta}_n, \zeta_0)[d(\tilde{\zeta}_n, \zeta_0) + \|\hat{\beta}_n - \beta_0\|]$; the second term on the right-hand side is bounded by $O_p(1)d^2(\tilde{\zeta}_n, \zeta_0)$; and the third term is bounded by $O_p(1)d(\tilde{\zeta}_n, \zeta_0)\|\hat{\beta}_n - \beta_0\|$. Thus, (10.27) follows from the fact that $d^2(\tilde{\zeta}_n, \zeta_0) = O_p(n^{-2/3})$ by Theorem 4.

7. $P_{\beta_0, h_0, H_0}\|\tilde{k}_{\hat{\beta}_n, \hat{h}_n, \widetilde{H}_n} - \tilde{k}_{\beta_0, h_0, H_0}\|^2 \to 0$ in probability as a consequence of the previously stated Lipschitz continuity of $\psi \mapsto \tilde{k}_\psi$ and consistency results in item 1 above. Furthermore, $P_{\hat{\beta}_n, h_0, H_0}\|\tilde{k}_{\hat{\beta}_n, \hat{h}_n, \widetilde{H}_n}\|^2 = O_p(1)$ from the boundedness assumptions. Thus, condition (10.28) of Theorem 7 below is satisfied.

Now all the conditions of Theorem 7 below are satisfied and, hence, $\hat{\beta}_n$ is efficient for $\beta_0$. $\square$

PROOF OF LEMMA 5. Fix $m > d$. Without loss of generality, we can assume by the i.i.d. structure that $X_i^* = X_i$ for $i = 1, \ldots, N_{m,n}$. Let $\epsilon_{0,n} \equiv \sqrt{n}(\tilde{\beta}_n - \beta_0) - n^{-1/2}\sum_{i=1}^n \phi_i$ and

$$
\epsilon_{j,n} \equiv (N_{m,n} - m)^{1/2}(\tilde{\beta}_{n,j}^* - \beta_0) - (N_{m,n} - m)^{-1/2}\sum_{i \in K_{j,n}} \phi_i,
$$

where $K_{j,n} \equiv \{1, \ldots, n\} - \{j, m + j, 2m + j, \ldots, (k_{m,n} - 1)m + j\}$, for $j = 1, \ldots, m$; and note that $\max_{0 \le k \le m}|\epsilon_{k,n}| = o_p(1)$ by asymptotic linearity. Now let $Z_{j,n}^* \equiv k_{m,n}^{-1/2}\sum_{i=1}^{k_{m,n}} \phi_{(i-1)m+j}$ for $j = 1, \ldots, m$, and define $\bar{Z}_n^* \equiv m^{-1}\sum_{j=1}^m Z_{j,n}^*$. Thus, $S_n^* = (m-1)^{-1}\sum_{j=1}^m (Z_{j,n}^* - \bar{Z}_n^*)(Z_{j,n}^* - \bar{Z}_n^*)^T + o_p(1)$. Hence, $S_n^*$ and $\sqrt{n}(\tilde{\beta}_n - \beta_0)$ are jointly asymptotically equivalent to $S_m \equiv$ $(m-1)^{-1} \sum_{j=1}^{m} (Z_j - \bar{Z}_m)(Z_j - \bar{Z}_m)^T$ and $Z_0$, respectively, where $\bar{Z}_m \equiv m^{-1} \sum_{j=1}^{m} Z_j$ and $Z_0, \ldots, Z_m$ are i.i.d. mean zero Gaussian deviates with variance $E[\phi\phi^T]$. Now the results follow by standard normal theory (see Appendix V of [39]). $\square$

LEMMA 8. *The class $\mathcal{F}_0$ in expression* (10.6) *is P-Donsker.*

PROOF. Let $k_0$ be the maximum possible value of $|\beta'Z + h(W)|$ whose existence is guaranteed by conditions A1(a), B2 and B3. Let $K$ in condition B5(c) be $[-k_0, k_0]$, and let $c_1$, $\alpha_1$ and $\xi_1$ be the choices of $c$, $\alpha$ and $\xi$ which satisfy the condition for this $K$. By condition B5(b), $F$ is one-to-one and, hence, $\mathcal{M}_{\mathbb{R}} = \{F^{-1}(\xi_1^{-1}G) : G \in \mathcal{M}_{[0,1]}\}$. Thus,

$$\mathcal{F}_0 = \left\{ \frac{\zeta(X; \beta, h, F^{-1}(\xi_1^{-1}G))}{1 + J(h)} : \beta \in \bar{B}_0, G \in \mathcal{M}_{[0,1]}, h \in \Im_\nu \right\},$$

where $\zeta$ is as defined in (10.3). Furthermore, for any $G_1, G_2 \in \mathcal{M}_{[0,1]}$, $\beta_1, \beta_2 \in \bar{B}_0$, and any $h_1, h_2 \in \Im_\nu$, we have

$$
\begin{aligned}
(10.26) \quad & \left| \frac{\zeta(X; \beta_1, h_1, F^{-1}(\xi_1^{-1}G_1))}{1 + J(h_1)} - \frac{\zeta(X; \beta_2, h_2, F^{-1}(\xi_1^{-1}G_2))}{1 + J(h_2)} \right| \\
& \leq m|\beta_1'Z - \beta_2'Z| \\
& \quad + \left| \frac{\zeta(X; \beta_1, h_1, F^{-1}(\xi_1^{-1}G_1)) - \zeta(X; \beta_1, h_1, F^{-1}(\xi_1^{-1}G_2))}{1 + J(h_1)} \right| \\
& \quad + \left| \frac{\zeta(X; \beta_2, h_1, F^{-1}(\xi_1^{-1}G_2))}{1 + J(h_1)} - \frac{\zeta(X; \beta_2, h_2, F^{-1}(\xi_1^{-1}G_2))}{1 + J(h_2)} \right|,
\end{aligned}
$$

where $m < \infty$ by B5(b) and the form of $\zeta$. By B5(c), the second term on the right-hand side of (10.26) is bounded above by $c_1|G_1(V) - G_2(V)|^{\alpha_1}$. Defining $\tilde{h}_j \equiv h_j / (1 + J(h_j))$, $j = 1, 2$, there exist constants $0 < c_2, c_3 < \infty$ so that, for the last term on the right-hand side of (10.26),

$$
\begin{aligned}
& \left| \frac{\zeta(X; \beta_2, h_1, F^{-1}(\xi_1^{-1}G_2))}{1 + J(h_1)} - \frac{\zeta(X; \beta_2, h_2, F^{-1}(\xi_1^{-1}G_2))}{1 + J(h_2)} \right| \\
& \leq \frac{c_2|h_1(W) - h_2(W)|}{(1 + J(h_1))} + c_3 \left| \frac{1}{1 + J(h_1)} - \frac{1}{1 + J(h_2)} \right| \\
& \leq c_2|\tilde{h}_1(W) - \tilde{h}_2(W)| + \{c_3 + c_2 h_2(W)\} \left| \frac{J(h_1) - J(h_2)}{(1 + J(h_1))(1 + J(h_2))} \right|,
\end{aligned}
$$

where $|h_2(W)| \leq c_0$ by constraint. Let $N_{[\cdot]}(\epsilon, \mathcal{F}, Q)$ be the bracketing number for the class $\mathcal{F}$ using $L_2(Q)$ brackets of size $\epsilon$. Note that the minimum number of points $S = \{s_1, \ldots, s_\nu\} \subset [0, \infty)$ needed to ensure that $\sup_{r \in [0,\infty)} \inf_{s \in S} |r - s| / [(1 + r)(1 + s)] < \epsilon$ is $O(\epsilon^{-1})$. Combining this with the facts that $\log N_{[\cdot]}(\epsilon, \mathcal{M}_{[0,1]}, Q) = c_4 \epsilon^{-1}$, where $c_4$ does not depend on $Q$ (Theorem 2.7.5 of [48]), and that both $\{h / \{1 + J(h)\} : h \in \Im_\nu\} \subset \mathcal{H} \equiv \{h \in \Im_\nu : J(h) \leq 1\}$ and $\log N_{[\cdot]}(\epsilon, \mathcal{H}, Q) = c_5 \epsilon^{-1}$, where $0 < c_5 < \infty$ does not depend on $Q$ (see Theorem 2.4 of [43]), we have that $\log N_{[\cdot]}(\epsilon, \mathcal{F}_0, P) \leq c_6 \epsilon^{-1/\alpha_1}$ for some $0 < c_6 < \infty$ and all $\epsilon \in (0, 1)$. Since this implies that the entropy integral with bracketing is bounded, the desired result follows. □

THEOREM 7 (Modification of Theorem 25.54 of [47]). *Suppose that the model $\{P_{\beta,\zeta} : \beta \in B_0\}$ is differentiable in quadratic mean with respect to $\beta$ at $(\beta_0, \zeta_0)$ and let the efficient information matrix $\tilde{I}_{\beta_0, \zeta_0}$ be nonsingular. Let $(\beta, \zeta) \mapsto \tilde{k}_{\beta, \zeta}$ be an estimating function satisfying $\tilde{k}_{\beta_0, \zeta_0} = \tilde{l}_{\beta_0, \zeta_0}$, where $\tilde{l}_{\beta_0, \zeta_0}$ is the efficient influence function for $\beta$ at $(\beta_0, \zeta_0)$. Let $\hat{\beta}_n$ satisfy $\sqrt{n} \mathbb{P}_n \tilde{k}_{\hat{\beta}_n, \hat{\zeta}_n} = o_p(1)$ and be consistent for $\beta_0$. In addition, suppose there exists a $P_{\beta_0, \zeta_0}$-Donsker class with square-integrable envelope function that contains every function $\tilde{k}_{\hat{\beta}_n, \hat{\zeta}_n}$ with probability tending to 1. Assume further that $\tilde{k}$ satisfies*

$$(10.27) \qquad \sqrt{n} P_{\hat{\beta}_n, \hat{\zeta}_0} \tilde{k}_{\hat{\beta}_n, \hat{\zeta}_n} = o_p(1 + \sqrt{n} \|\hat{\beta}_n - \beta_0\|)$$

*and*

$$(10.28) \quad P_{\beta_0, \zeta_0} \|\tilde{k}_{\hat{\beta}_n, \hat{\zeta}_n} - \tilde{k}_{\beta_0, \zeta_0}\|^2 = o_p(1), \qquad P_{\hat{\beta}_n, \hat{\zeta}_0} \|\tilde{k}_{\hat{\beta}_n, \hat{\zeta}_n}\|^2 = O_p(1).$$

*Then $\hat{\beta}_n$ is asymptotically efficient at $(\beta_0, \zeta_0)$.*

PROOF. The proof is almost identical to van der Vaart's proof of his Theorem 25.54 in [47]. □

**Acknowledgments.** The authors thank Drs. Yi Lin and Jason P. Fine for several insightful discussions. The authors also express appreciation to the referees and editors for several helpful suggestions.

Department of Statistics
University of Wisconsin–Madison
1300 University Avenue
Madison, Wisconsin 53706
USA
e-mail: shuangge@biostat.wisc.edu

Departments of Statistics and
    Biostatistics and Medical Informatics
University of Wisconsin–Madison
1300 University Avenue
Madison, Wisconsin 53706
USA
e-mail: kosorok@biostat.wisc.edu